\documentclass{amsart}
\usepackage{graphicx}

\newtheorem{theorem}{Theorem}[section]

\newtheorem{corollary}[theorem]{Corollary}
\newtheorem{proposition}[theorem]{Proposition}

\theoremstyle{definition}
\newtheorem{definition}[theorem]{Definition}

\newtheorem{conjecture}[theorem]{Conjecture}

\theoremstyle{remark}

\numberwithin{equation}{section}

\begin{document}

\title{Ascending number of knots and links}

\author{Makoto Ozawa}
\address{Department of Natural Sciences, Faculty of Arts and Sciences, Komazawa University, 1-23-1 Komazawa, Setagaya-ku, Tokyo, 154-8525, Japan}
\email{w3c@komazawa-u.ac.jp}





\begin{abstract}
We introduce a new numerical invariant of knots and links from the descending diagrams.
It is considered to live between the unknotting number and the bridge number.
\end{abstract}

\maketitle

\section{Introduction}
Throughout this paper we work in the piecewise linear category.
We shall study knots and links in the three-dimensional Euclidean space $\Bbb{R}^{3}$.
For the standard definitions and results of knots and links, we refer to \cite{A}, \cite{BZ}, \cite{C}, \cite{Kaw}, \cite{Lic2}, \cite{L}, \cite{M} and \cite{R}.

For the purpose of defining a numerical invariant of knots and links, we need to prepare following terms.
A link diagram is {\em ordered} if an order is given to its components.
A link diagram is {\em based} if a basepoint (different from the crossing points) is specified on each component.
A link diagram is {\em oriented} if an orientation is specified on each component.

\begin{figure}[htbp]
	\begin{center}
	\begin{tabular}{cc}
	\includegraphics[trim=0mm 0mm 0mm 0mm, width=.4\linewidth]{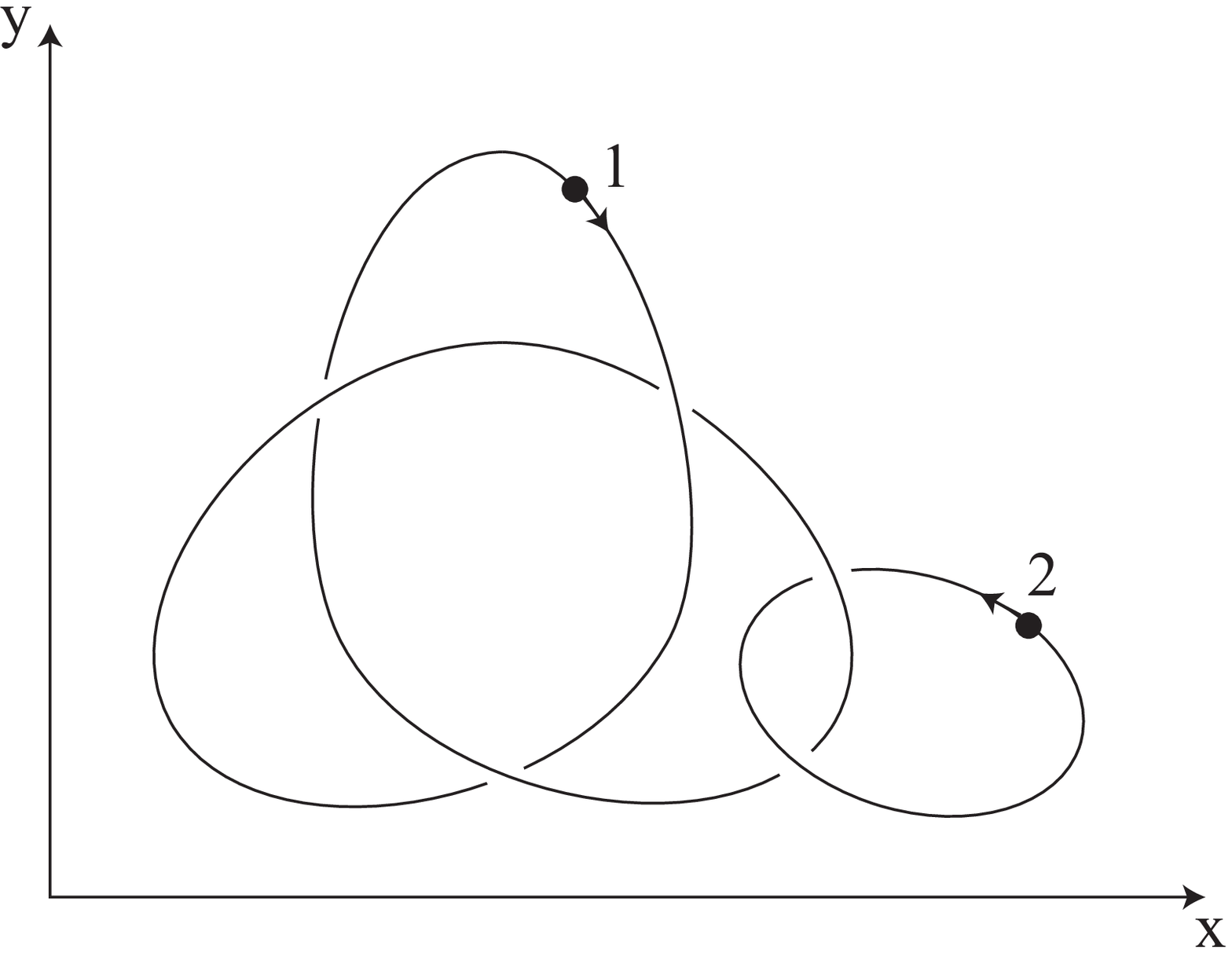}&
	\includegraphics[trim=0mm 0mm 0mm 0mm, width=.4\linewidth]{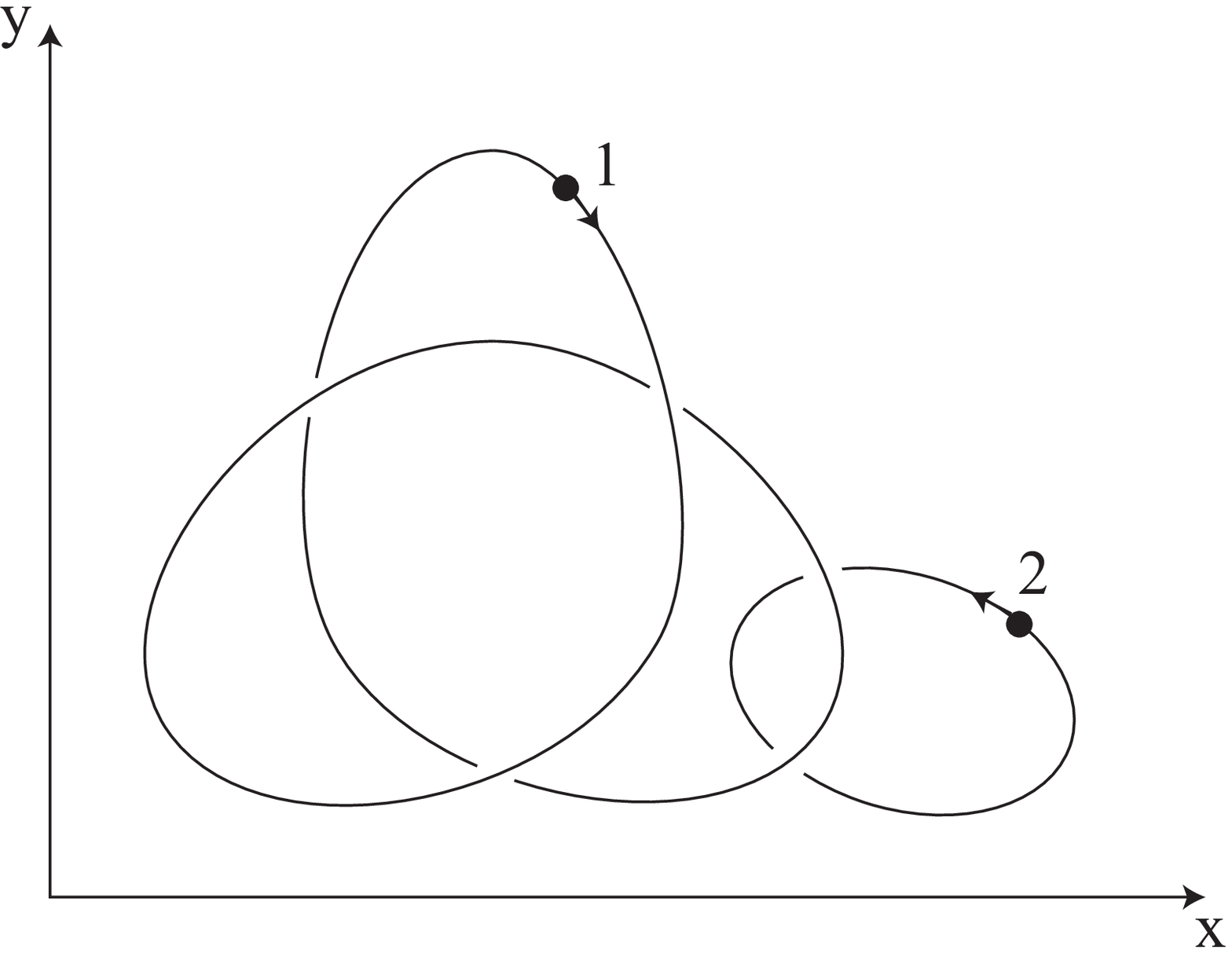}\\
	$\tilde{L}$ & $d(\tilde{L})$
	\end{tabular}
	\end{center}
	\caption{}
	\label{example1}
\end{figure}

Let $L$ be a link, and let $\tilde{L}$ be a based ordered oriented link diagram of $L$.
The {\it descending diagram} of $\tilde{L}$ is obtained as follows, and will be denoted by $d(\tilde{L})$.
Beginning at the basepoint of the first component of $\tilde{L}$ and proceeding in the direction specified by the orientation, change the crossings as necessary so that each crossing is first encountered as an over-crossing.
Continue this procedure with the remaining components in the sequence determined by the ordering, proceeding from the basepoint in the direction determined by the orientation, changing crossings so that ultimately every crossing is first encountered as an over-crossing.
The result is the descending diagram $d(\tilde{L})$ obtained from $\tilde{L}$. An example is shown in Figure \ref{example1}.
Note that $d(\tilde{L})$ is a diagram of a trivial link.

Prof. Shin'ichi Suzuki mentioned in \cite[2.16 (ii)]{S} that the smallest number of crossing changes to obtain the trivial link may not always be provided with the number of difference crossing between a link diagram and the descending diagram.
That was why we defined the next numerical invariant in \cite{MO}.

\begin{definition}
Let $L$ be a link and let $\tilde{L}$ be a based ordered oriented link diagram of $L$. The {\em ascending number} of $\tilde{L}$ is defined as the number of different crossings between $\tilde{L}$ and $d(\tilde{L})$, and denoted by $a(\tilde{L})$. 
The {\em ascending number} of $L$ is defined as the minimum number of $a(\tilde{L})$ over all based ordered oriented link diagram $\tilde{L}$ of $L$, and denoted by $a(L)$.
\end{definition}

The definition and the above note imply that $a(L) \geq u(L)$,
where $u(L)$ denotes the unknotting number of $L$. 
Accordingly $L$ is a trivial link if and only if $a(L) = 0$.
A difference of $a(L)$ and $u(L)$ can become no matter how big.
Actually, the mention of S. Suzuki will be demonstrated in Corollary \ref{difference}.

Now this paper is organized as follows. In section 2 we
state all results. In section 3 we give all proofs of the
results. In section 4 we show a table of the ascending
number, the unknotting number, and the bridge number of prime knots with eight crossings or less.

\section{Results}

\begin{proposition}\label{crossing}
For a nontrivial knot $K$, we have
\[ a(K) \leq \left[\frac{c(K)-1}{2}\right].\]
For a link $L$, we have 
\[ a(L) \leq \left[\frac{c(L)}{2}\right],\] 
where $c(L)$ denotes the minimum crossing number of $L$,
and $[x]$ denotes the greatest integer which does not
exceed $x$.
\end{proposition}
%

The following theorem is the fundamental inequality between the ascending
number and the bridge number.

\begin{theorem}\label{bridge}
For an $n$-component link $L$, we have
\[a(L) \geq b(L)-n,\]
where $b(L)$ is the bridge number of $L$.
\end{theorem}

The following corollary asserts the difference between
$a(K)$ and $u(K)$.

\begin{corollary}\label{difference}
For any nonnegative integer $n$, there is a knot $K$
such that $a(K)-u(K) \geq n$.
\end{corollary}

For the connected sum, we have the following proposition. 

\begin{proposition}\label{subadditive}
For a knot, the ascending number is subadditive with
respect to connected sum \#, i.e. $a(K_{1}\#K_{2}) \leq
a(K_{1})+a(K_{2}).$
\end{proposition}

It is natural to ask whether the ascending number is
additive with respect to the connected sum. 

\begin{conjecture}
For any knots $K_{1}, K_{2}$,
\[a(K_{1}\# K_{2}) = a(K_{1}) + a(K_{2}).\]
\end{conjecture}

The following corollary solves the above question
partially.

\begin{corollary}\label{subanswer}
Let $K_{i}\ (i = 1,2)$ be a knot. If $K_{i}$ satisfies
$a(K_{i}) = b(K_{i}) - 1\ (i = 1,2)$, then
\[a(K_{1}\# K_{2}) = a(K_{1}) + a(K_{2}).\]
\end{corollary}

In consequence of this, we have the following corollary. 

\begin{corollary}\label{exist}
For any nonnegative integer $n$, there is a knot $K$ such that
$a(K) = n$.
\end{corollary}

The following theorem characterizes the ascending number one link.
We note that Tat Sang Fung (\cite{Fun}) also obtained the same result for knot case.
In fact, he showed that almost descending knots, i.e. knots with diagrams descending except at one crossing, are twist knots.

\begin{theorem}\label{one}
Let $L$ be an $n$-component link. Then the ascending
number of $L$ is one if and only if $L$ is 
\begin{align*}
& (\text{twist knot})\circ O^{n-1},
 \text{\quad if L is completely splittable}\\
\intertext{or}
& (\text{Hopf link})\circ O^{n-2},
 \text{\quad otherwise,}
\end{align*}
where $\circ$ denotes the split union, and $O^{n}$
denotes an $n$-component trivial link.

In paticular, the ascending number of a knot $K$ is one if
and only if $K$ is a twist knot.
\end{theorem}

The following theorem determines the ascending number for
torus knots.

\begin{theorem}\label{torus}

Let $p$ and $q$ be coprime positive integers, and let
$T_{p,q}$ be a $(p,q)$-torus knot. Then we have
\[ a(T_{p,q})=\frac{(p-1)(q-1)}{2}. \] 
\end{theorem}




\section{Proofs}

\begin{proof}(of Proposition \ref{crossing})
At the beginning, we prove the proposition for a
link. Let $L$ be a link, and let $\tilde{L}$ be a based
ordered oriented link diagram of $L$ with minimum crossings.
Consider a based ordered oriented link diagram which is
obtained by reversing the order and the orientation of
$\tilde{L}$, and let it be denoted by $-\tilde{L}$. Then
it holds that $c(L) = a(\tilde{L}) + a(-\tilde{L})$.
Indeed $d(-\tilde{L})$ is obtained by reflecting
$d(\tilde{L})$ in the plane. Hence
\[\frac{c(L)}{2} \geq \min\{a(\tilde{L}),a(-\tilde{L})\} \geq
a(L).\]

Next we prove the proposition for a knot. Let $K$ be a
knot, and let $\tilde{K}$ be a knot diagram of $K$ with
minimum crossings. Choose a crossing of $\tilde{K}$ and
specify a basepoint and an orientation on $\tilde{K}$ in
the following fashion. 
Whenever we begin the basepoint of $\tilde{K}$ and proceed
in the direction specified by the orientation, we
first encounter the crossing as an over-crossing.  

On the other hand let $-\tilde{K'}$ be a based oriented
knot diagram which is obtained by reversing the
orientation of $\tilde{K}$ and sliding the base point of
$\tilde{K}$, as follows. Whenever we begin the basepoint of
$-\tilde{K'}$ and proceed in the direction specified by the
orientation, we first encounter the crossing as an
over-crossing. See Figure \ref{reverse}.

\begin{figure}[htbp]
	\begin{center}
	\begin{tabular}{cc}
	\includegraphics[trim=0mm 0mm 0mm 0mm, width=.25\linewidth]{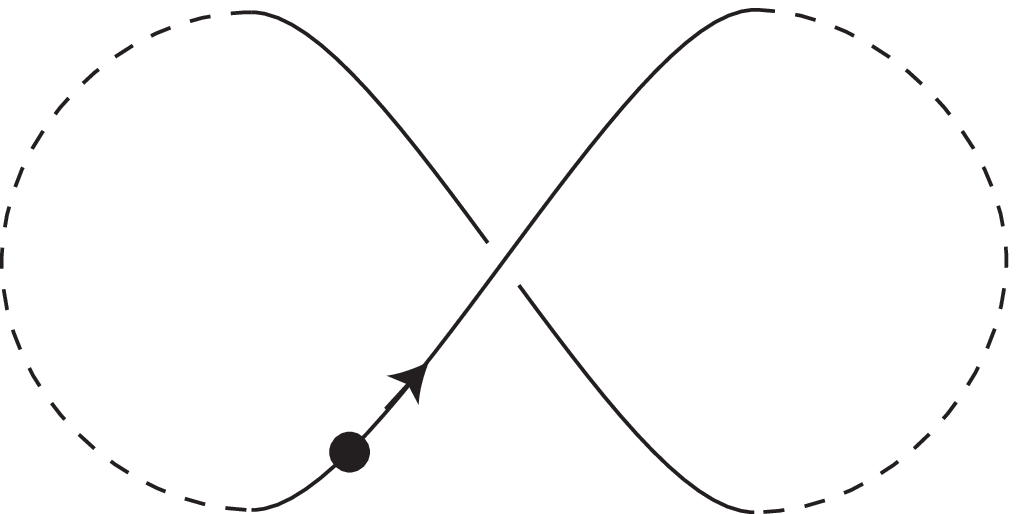}&
	\includegraphics[trim=0mm 0mm 0mm 0mm, width=.25\linewidth]{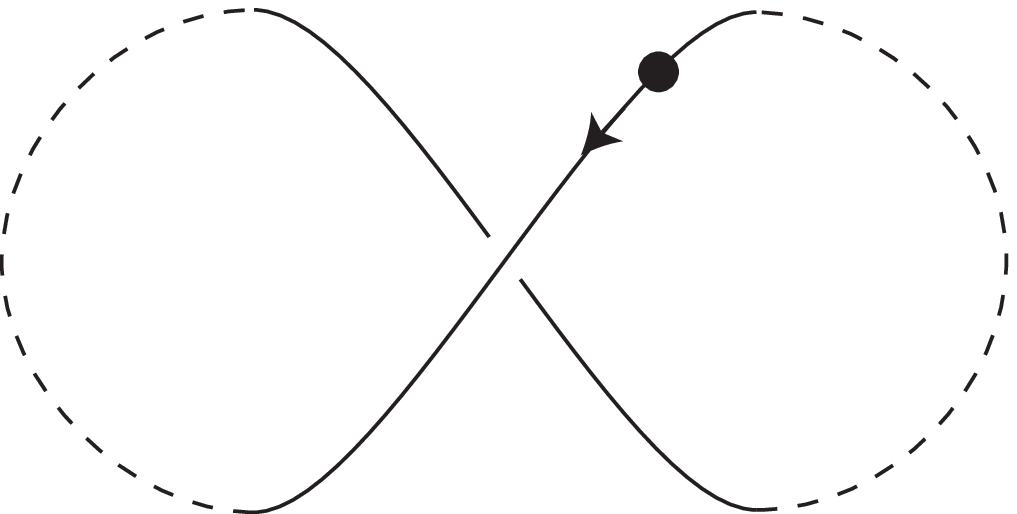}\\
	$\tilde{K}$ & $-\tilde{K}'$
	\end{tabular}
	\end{center}
	\caption{}
	\label{reverse}
\end{figure}

Then it holds that $c(K) - 1 = a(\tilde{K}) +
a(-\tilde{K'})$. Indeed a knot consists of a single circle.
Hence 

\[
\frac{c(k)-1}{2} \geq \min\{a(\tilde{K}),a(-\tilde{K'})\}
\geq a(K).
\]

This completes the proof.
\end{proof}

\begin{proof}(of Theorem \ref{bridge})
Let $\tilde{L}$ be a based ordered oriented link diagram of $L=K_1\cup \cdots \cup K_n$ satisfying $a(\tilde{L})=a(L)$.

First, we make an $n$-bridge presentation of a trivial link $L'=K_1'\cup \cdots \cup K_n'$ from the descending diagram $d(\tilde{L})$ by isotoping on $z$-coordinate so that $L'$ satisfy the next conditions.
We note that a similar presentation can be seen in \cite[Section 3.1]{A}.

\begin{enumerate}
\item For each $i\in \{1,\ldots,n\}$, $K_i'$ is parametrized by an embedding $f_i:[0,1]\to \Bbb{R}^3$ so that the orientation and basepoint of $p(f_i([0,1]))$ are in consistency with that of $d(\tilde{K_i})$, where $p$ denotes the projection $\Bbb{R}^3\to \Bbb{R}^2\times \{0\},\ (x,y,z)\mapsto (x,y,0)$.
\item For any $x,\ y\in [0,1]$ such that $x<y$, $h(f_i(x))>h(f_i(y))$, where $h$ denotes the hight function $\Bbb{R}^3\to \{0\}^2\times \Bbb{R},\ (x,y,z)\mapsto (0,0,z)$.
\item For each $i\in \{1,\ldots,n-1\}$, $h(f_i(\{1\}))>h(f_{i+1}(\{0\}))$.
\item For each $i\in \{1,\ldots,n\}$, $f_i(\{0\})$ and $f_i(\{1\})$ are connected by a vertical arc $P_i$ contained in $p^{-1}($base point of $d(\tilde{K_i}))$. Thus, $K_i'=f_i([0,1])\cup P_i$
\end{enumerate}

Next, we deform $L'$ to obtain an $(n+a(L))$-bridge presentation of $L$ as follows.
Let $\mathcal{C}$ be the set of crossings of $\tilde{L}$ which are different from that of $d(\tilde{L})$.
For each crossing $c\in \mathcal{C}$, we isotope the corresponding under-crossing of $L'$ to exceed the corresponding over-crossing.
Here, we may assume that such an operation accrues exactly one maximal point.
As a result, we obtain a presentation of $L$ whose bridge number is equal to $n+a(L)$.
Figure \ref{example2} shows a realization of $\tilde{L}$ and $d(\tilde{L})$ in Figure \ref{example1}.

\begin{figure}[htbp]
	\begin{center}
	\begin{tabular}{cc}
	\includegraphics[trim=0mm 0mm 0mm 0mm, width=.3\linewidth]{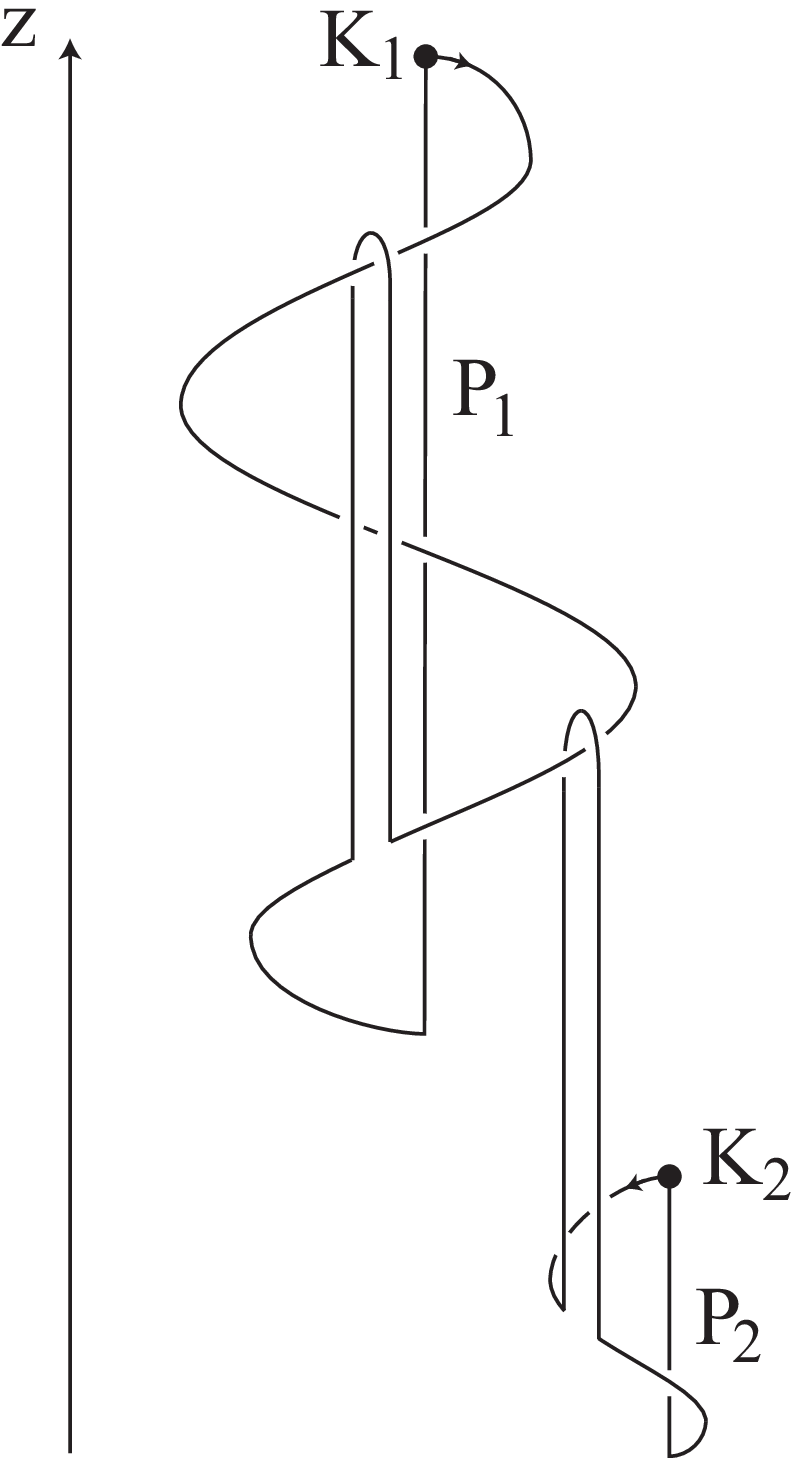}&
	\includegraphics[trim=0mm 0mm 0mm 0mm, width=.3\linewidth]{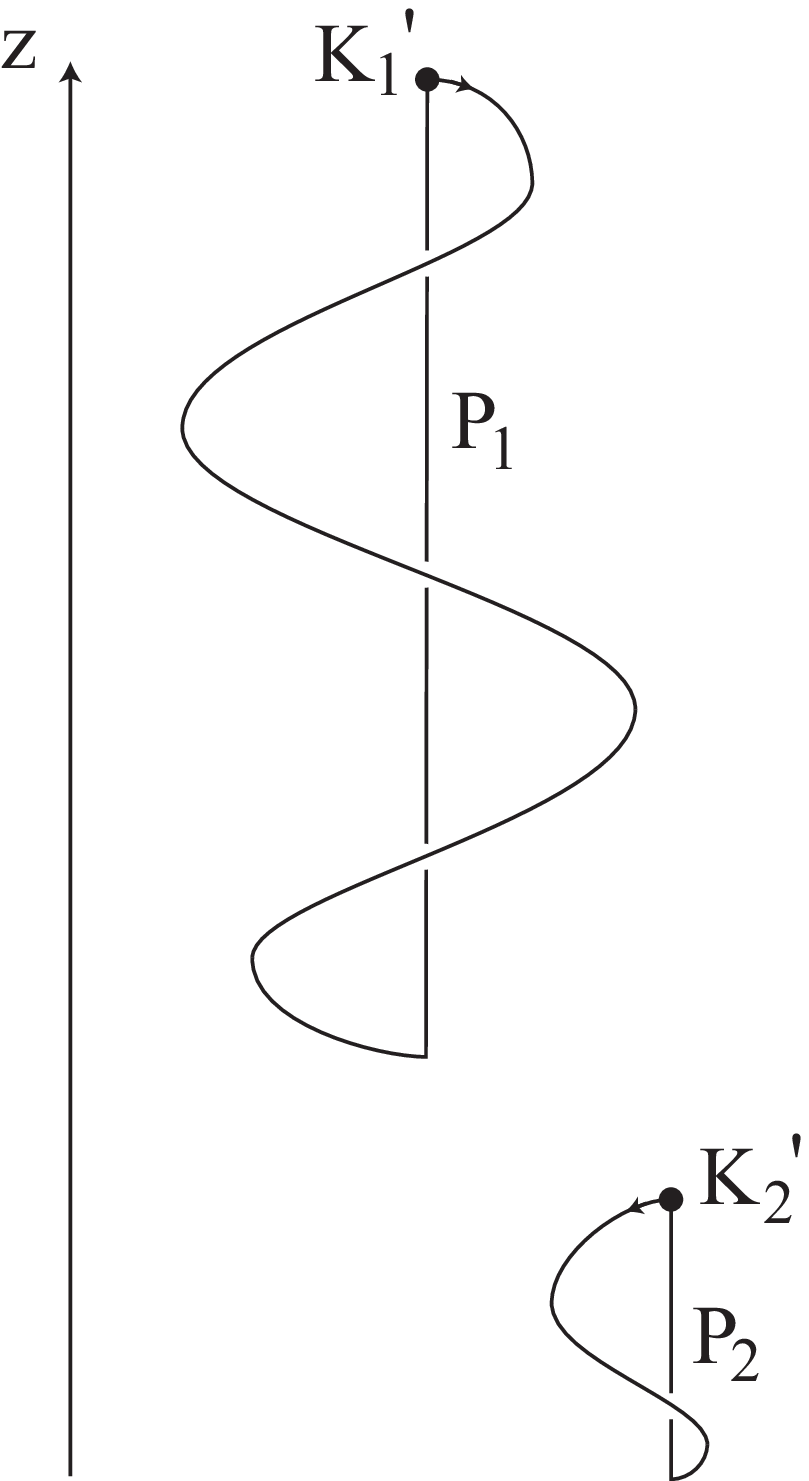}\\
	$L$ & $L'$
	\end{tabular}
	\end{center}
	\caption{}
	\label{example2}
\end{figure}

Hence we have
\[
b(L)\le n+a(L).
\]

This arrives at the inequality of Theorem \ref{bridge}.
\end{proof}

\begin{proof}(of Corollary \ref{difference})
We will show that a composite knot $K$ which is a connected sum of $n$ $8_{17}$'s satisfies the inequality of Corollary \ref{difference}.
It is known that $u(8_{17})=1$ and $b(8_{17})=3$.
By Schubert's formula $b(K_1\#K_2)=b(K_1)+b(K_2)-1$ (\cite{Sch}),
\[
b(K)=nb(8_{17})-(n-1)=2n+1.
\]

Since an inequality $u(K_1\#K_2)\le u(K_1)+u(K_2)$ generally holds,
\[
u(K)\le n u(8_{17})=n.
\]

Hence, by Theorem \ref{bridge},
\[
a(K)-u(K)\ge (b(K)-1)-u(K)\ge 2n-n=n.
\]
\end{proof}

\begin{proof}(of Proposition \ref{subadditive})
Let $\tilde{K_1}$ and $\tilde{K_2}$ be based oriented diagram of knots $K_1$ and $K_2$ such that $a(\tilde{K_i})=a(K_i)$ for $i=1,2$.
By considering diagrams on the 2-sphere, we may assume that basepoints of $K_1$ and $K_2$ are in regions containing a point at infinity.
Then we perform a connected sum $\tilde{K_1}$ and $\tilde{K_2}$ to obtain a based oriented diagram $\tilde{K_1\# K_2}$ as Figure \ref{connected sum}.

\begin{figure}[htbp]
	\begin{center}
	\includegraphics[trim=0mm 0mm 0mm 0mm, width=.7\linewidth]{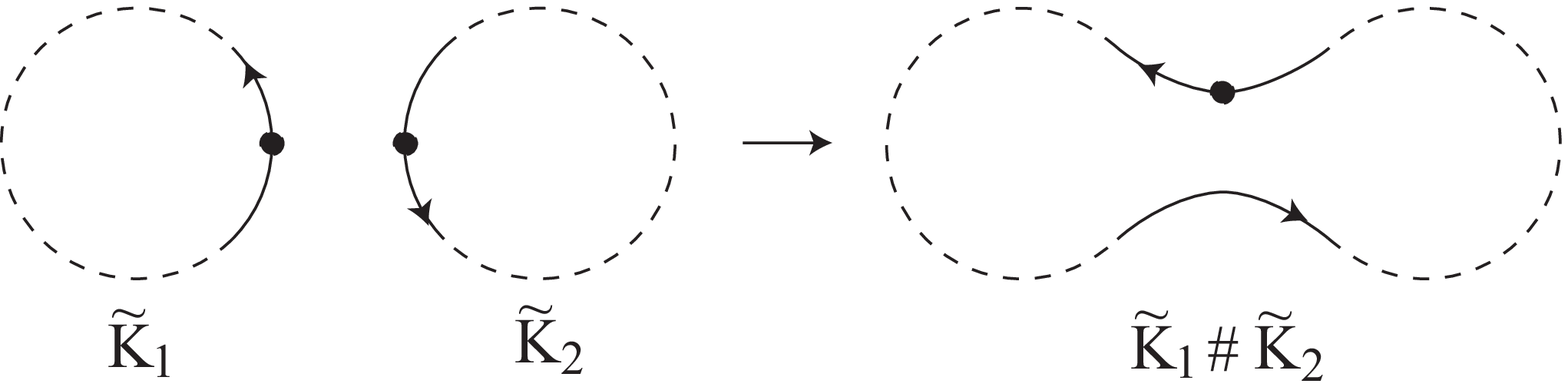}	\end{center}
	\caption{}
	\label{connected sum}
\end{figure}

Therefore
\[
a(K_1\# K_2)\le a(\tilde{K_1}\#\tilde{K_2})=a(\tilde{K_1})+a(\tilde{K_2})=a(K_1)+a(K_2).
\]
\end{proof}

\begin{proof}(of Corollary \ref{subanswer})
By Schubert's formula and Theorem \ref{bridge},
\[
a(K_1\# K_2)\ge b(K_1\# K_2)-1 =(b(K_1)-1)+(b(K_2)-1)=a(K_1)+a(K_2).
\]
Then, Proposition \ref{subadditive} gives the equality.
\end{proof}

\begin{proof}(of Corollary \ref{exist})
Let $K$ be a connected sum of $n$ $3_1$'s.
Since $a(3_1)=1$ and $b(3_1)=2$, $3_1$ satisfies the supposition of Corollary \ref{subanswer}.
Hence, we have $a(3_1)=n$.
\end{proof}

\begin{proof}(of Theorem \ref{one})
The ascending number of the Hopf link is obviously equal to 1.
The ascending number of a twist knot is also equal to 1 by the transformation in Figure \ref{twist}.

\begin{figure}[htbp]
	\begin{center}
	\includegraphics[trim=0mm 0mm 0mm 0mm, width=.8\linewidth]{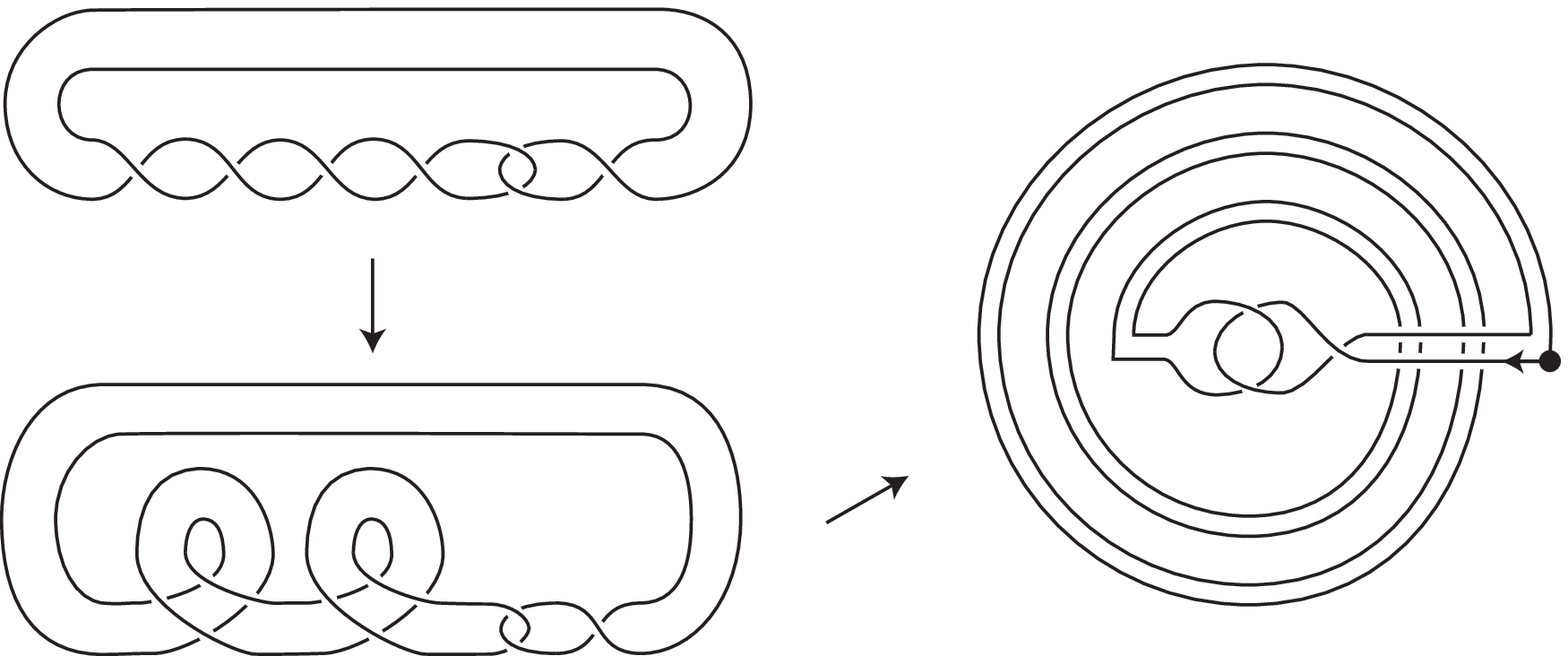}	\end{center}
	\caption{}
	\label{twist}
\end{figure}

Conversely, suppose that $a(L)=1$ and let $\tilde{L}$ be a based ordered oriented diagram of $L$ with $a(\tilde{L})=1$.
Then there is just one crossing $c$ of $\tilde{L}$ that is a difference between $\tilde{L}$ and $d(\tilde{L})$.
We make an $n$-bridge presentation of a trivial link $L'$ from the descending diagram $d(\tilde{L})$ and deform $L'$ to obtain an $(n+1)$-bridge presentation of $L$ in the same way as the proof of Theorem \ref{bridge}.

There are two cases.

\begin{enumerate}
\item The crossing $c$ consists of an over-crossing and an under-crossing of $\tilde{K_i}$.
\item The crossing $c$ consists of an over-crossing of $\tilde{K_j}$ and an under-crossing of $\tilde{K_i}$ $(i<j)$.
\end{enumerate}

In Case 1, $L$ is completely splittable and by a deformation in Figure \ref{twist2}, $K_i$ is a twist knot.
Here, we may assume that a basepoint of $\tilde{K_i}$ is in the left most by considering a diagram on the 2-sphere.

\begin{figure}[htbp]
	\begin{center}
	\includegraphics[trim=0mm 0mm 0mm 0mm, width=.9\linewidth]{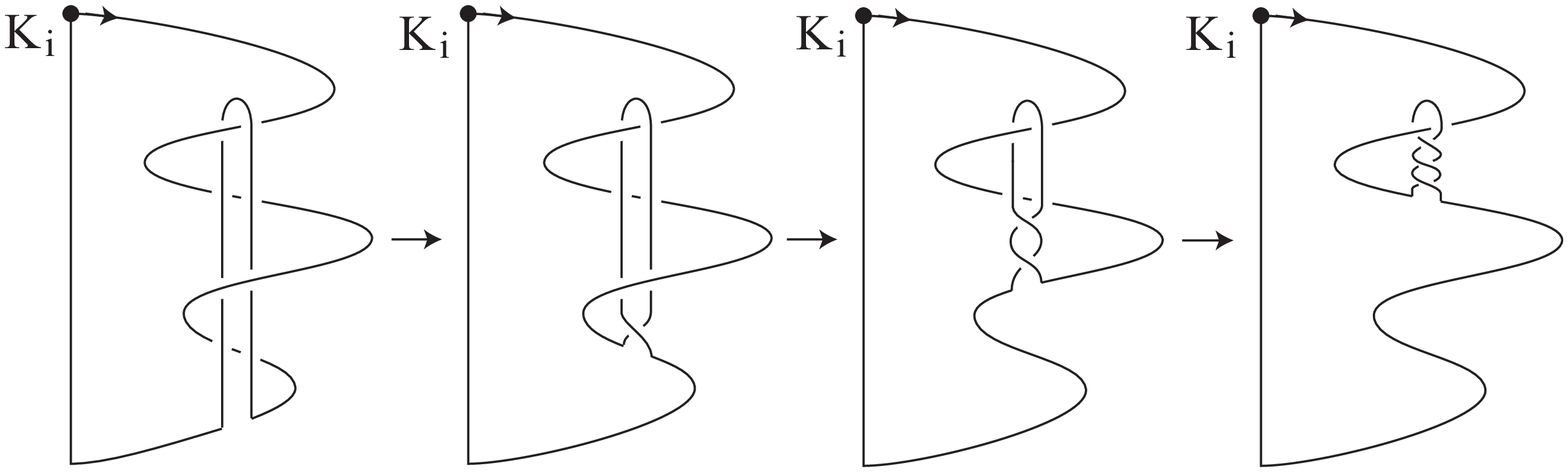}	\end{center}
	\caption{}
	\label{twist2}
\end{figure}

In Case 2, $L$ is a split union of a Hopf link $K_i\cup K_j$ and $(n-2)$-component trivial link by a deformation in Figure \ref{Hopf2}.

\begin{figure}[htbp]
	\begin{center}
	\includegraphics[trim=0mm 0mm 0mm 0mm, width=.9\linewidth]{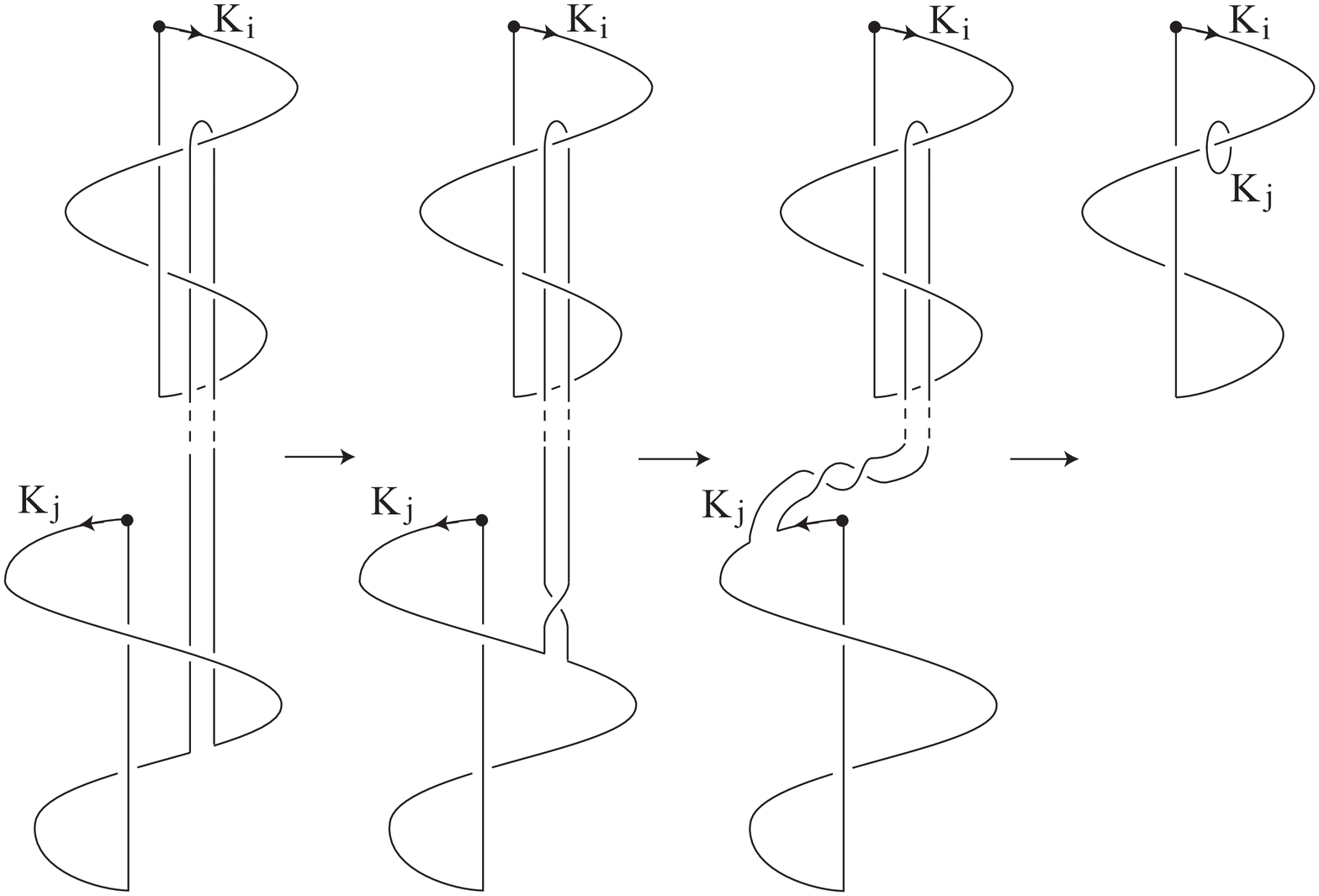}	\end{center}
	\caption{}
	\label{Hopf2}
\end{figure}

This completes the proof of Theorem \ref{one}.
\end{proof}

\begin{proof}(of Theorem \ref{torus})
Let $p,\ q$ be coprime positive integers $(p<q)$, $K$ a $(p,q)$-torus knot and $\tilde{K}$ a standard diagram with $(p-1)q$-crossings.

We assign a basepoint and an orientation to $\tilde{K}$ so that whenever we begin the basepoint and proceed by the orientation, we first encounter $(p-1)$-crossings as over-crossings.
We denote a resultant based oriented diagram by $\tilde{K_a}$ (Figure \ref{torus2}).

On the other hand, let $\tilde{K_b}$ be a based oriented diagram which is obtained by reversing the orientation of $\tilde{K_a}$ and sliding the base point of $\tilde{K_a}$ so that whenever we begin the basepoint of $\tilde{K_b}$ and proceed by the orientation, we first encounter the same $(p-1)$-crossings as over-crossings (Figure \ref{torus2}).

\begin{figure}[htbp]
	\begin{center}
	\includegraphics[trim=0mm 0mm 0mm 0mm, width=.5\linewidth]{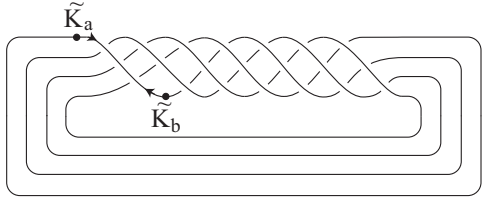}	\end{center}
	\caption{$(4,7)$-torus knot}
	\label{torus2}
\end{figure}

Then, it holds that $a(\tilde{K_a})+a(\tilde{K_b})=(p-1)(q-1)$.
Hence
\[
a(\tilde{K})\le \min \{ a(\tilde{K_a}), a(\tilde{K_b}) \} \le \frac{(p-1)(q-1)}{2}
\]

On the other hand, by a solution of Milnor's conjecture $u(K)=(p-1)(q-1)/2$ (\cite{KM}),
\[
\frac{(p-1)(q-1)}{2} = u(K) \le a(K) \le a(\tilde{K})
\]
Hence, we obtain Theorem \ref{torus}.
\end{proof}

\section{Table of the ascending number of knots}

Here we present a table of the ascending number, unknotting number and bridge number of knots with $8$-crossings or less.
The unknotting numbers are cited from \cite{S} except for $u(8_{10})=2$ from \cite[Corollary 1.2]{OS} and $u(8_{16})=2$ from an unpublished result by J. R. Rickard (\cite{Lic}, \cite{K}).
The bridge numbers are cited from \cite{M}.
For the ascending numbers, $a(8_4)=2$ was determined by Miki Okuda (\cite{O}), $a(8_{13})=2$ was determined by Sachie Fujimura (\cite{F}) and she also pointed out $a(8_{20})=2$, where they belong to Laboratory of Prof. Takao Matumoto, Department of Mathematics, School of Science, Hiroshima University.
Miki Okuda also determined $a(9_3)=3$, $a(9_4)=2$, $a(9_6)=3$, $a(9_7)=2$, $a(9_{47})=2$ and $a(9_{48})=2$.

\bigskip
\begin{minipage}{4cm}
\begin{tabular}{|c||c|c|c|}
\hline
$K$ & $a(K)$ & $u(K)$ & $b(K)$\\
\hline
$3_1$ & 1 & 1 & 2\\
\hline
$4_1$ & 1 & 1 & 2\\
\hline
$5_1$ & 2 & 2 & 2\\
\hline
$5_2$ & 1 & 1 & 2\\
\hline
$6_1$ & 1 & 1 & 2\\
\hline
$6_2$ & 2 & 1 & 2\\
\hline
$6_3$ & 2 & 1 & 2\\
\hline
\end{tabular}
\end{minipage}
\hfill
\begin{minipage}{6cm}
\begin{tabular}{|c||c|c|c|}
\hline
$K$ & $a(K)$ & $u(K)$ & $b(K)$\\
\hline
$7_1$ & 3 & 3 & 2\\
\hline
$7_2$ & 1 & 1 & 2\\
\hline
$7_3$ & 2 & 2 & 2\\
\hline
$7_4$ & 2 & 2 & 2\\
\hline
$7_5$ & 2 & 2 & 2\\
\hline
$7_6$ & 2 & 1 & 2\\
\hline
$7_7$ & 2 & 1 & 2\\
\hline
\end{tabular}
\end{minipage}
\bigskip

\begin{minipage}{4cm}
\begin{tabular}{|c||c|c|c|}
\hline
$K$ & $a(K)$ & $u(K)$ & $b(K)$\\
\hline
$8_1$ & 1 & 1 & 2\\
\hline
$8_2$ & 2 or 3 & 2 & 2\\
\hline
$8_3$ & 2 & 2 & 2\\
\hline
$8_4$ & 2 & 2 & 2\\
\hline
$8_5$ & 2 or 3 & 2 & 3\\
\hline
$8_6$ & 2 & 2 & 2\\
\hline
$8_7$ & 2 or 3 & 1 & 2\\
\hline
$8_8$ & 2 & 2 & 2\\
\hline
$8_9$ & 2 or 3 & 1 & 2\\
\hline
$8_{10}$ & 2 or 3 & 2 & 3\\
\hline
\end{tabular}
\end{minipage}
\hfill
\begin{minipage}{6cm}
\begin{tabular}{|c||c|c|c|}
\hline
$K$ & $a(K)$ & $u(K)$ & $b(K)$\\
\hline
$8_{11}$ & 2 & 1 & 2\\
\hline
$8_{12}$ & 2 & 2 & 2\\
\hline
$8_{13}$ & 2 & 1 & 2\\
\hline
$8_{14}$ & 2 & 1 & 2\\
\hline
$8_{15}$ & 2 & 2 & 3\\
\hline
$8_{16}$ & 2 or 3 & 2 & 3\\
\hline
$8_{17}$ & 2 or 3 & 1 & 3\\
\hline
$8_{18}$ & 2 & 2 & 2\\
\hline
$8_{19}$ & 3 & 3 & 3\\
\hline
$8_{20}$ & 2 & 1 & 3\\
\hline
$8_{21}$ & 2 & 1 & 3\\
\hline
\end{tabular}
\end{minipage}
\bigskip

We present a based oriented knot diagram which gives the minimal ascending number.
We omit knot diagrams which are same as minimal crossing diagram.
As S. Suzuki remarked in \cite[2.16 (i)]{S}, the minimal ascending diagram is not unique.

\begin{figure}[htbp]
	\begin{center}
	\begin{tabular}{ccc}
	\includegraphics[trim=0mm 0mm 0mm 0mm, width=.28\linewidth]{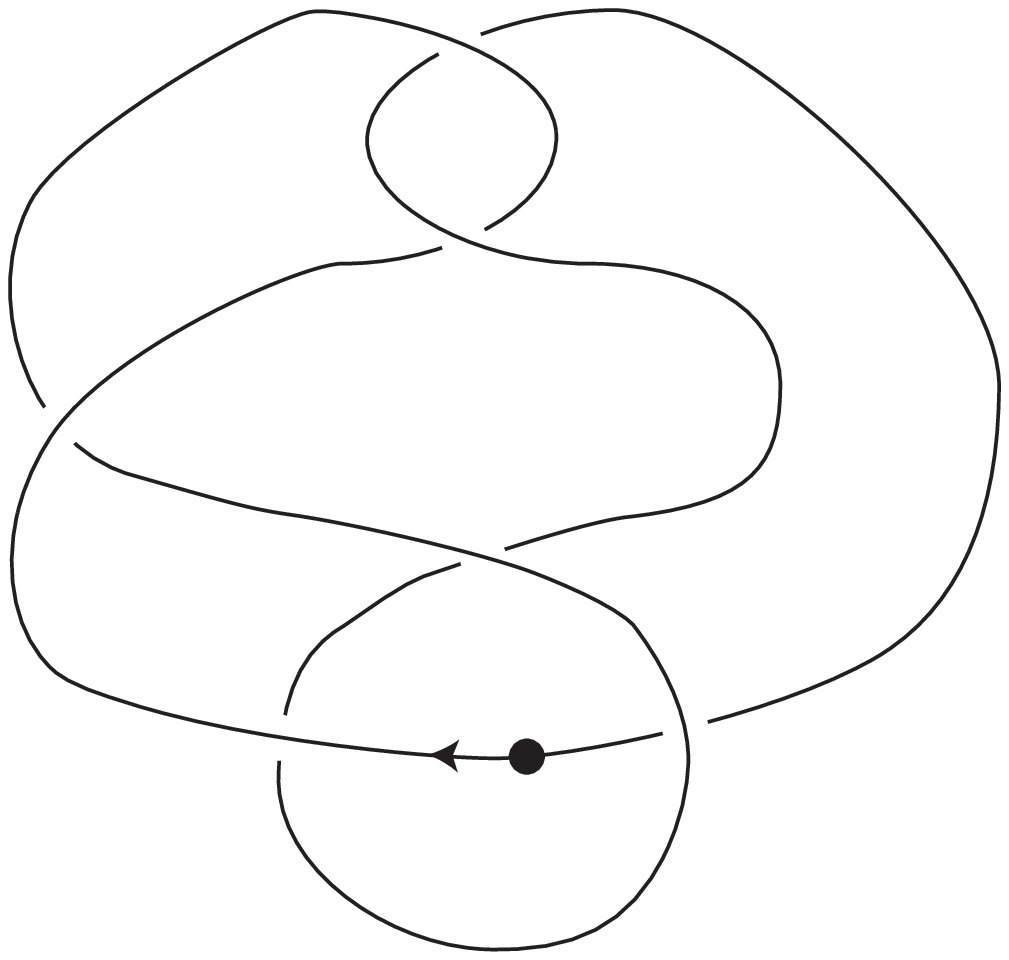}&
	\includegraphics[trim=0mm 0mm 0mm 0mm, width=.3\linewidth]{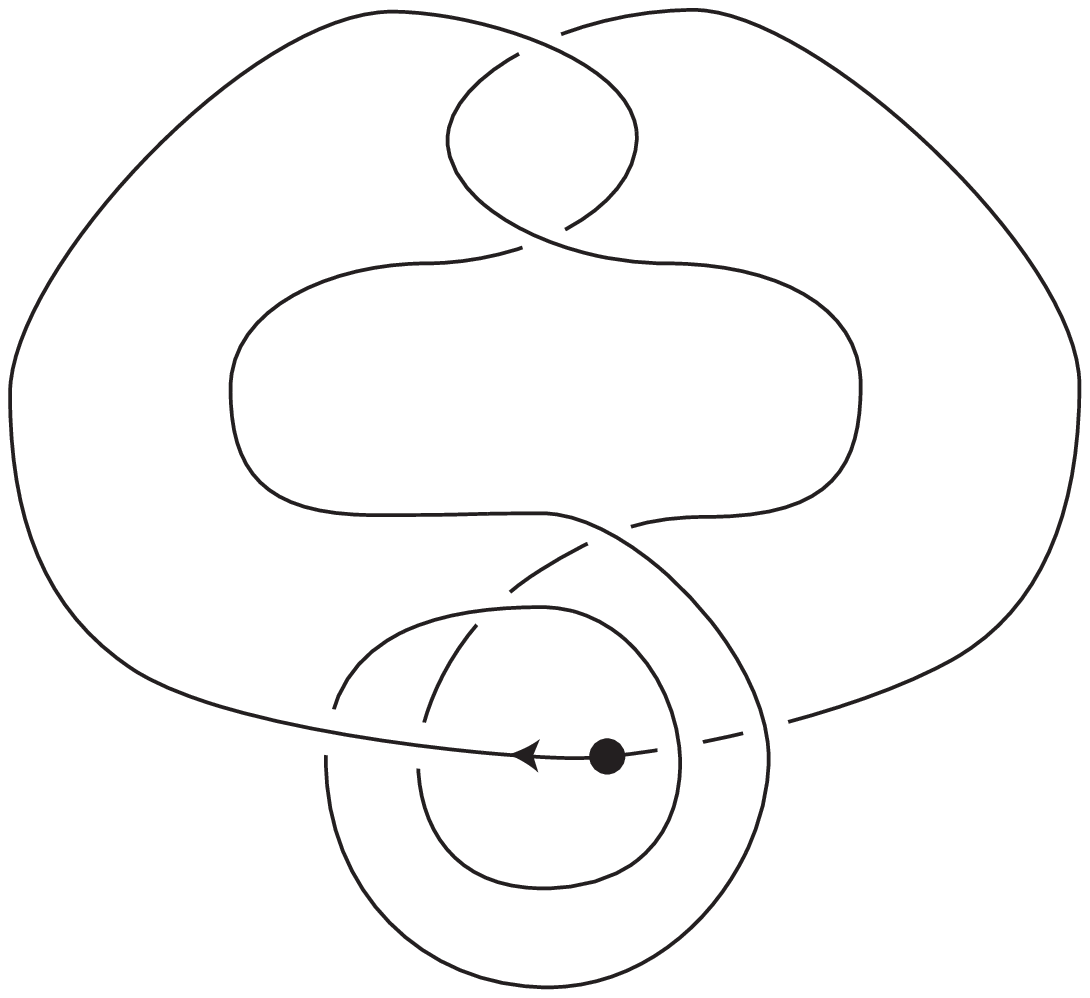}&
	\includegraphics[trim=0mm 0mm 0mm 0mm, width=.29\linewidth]{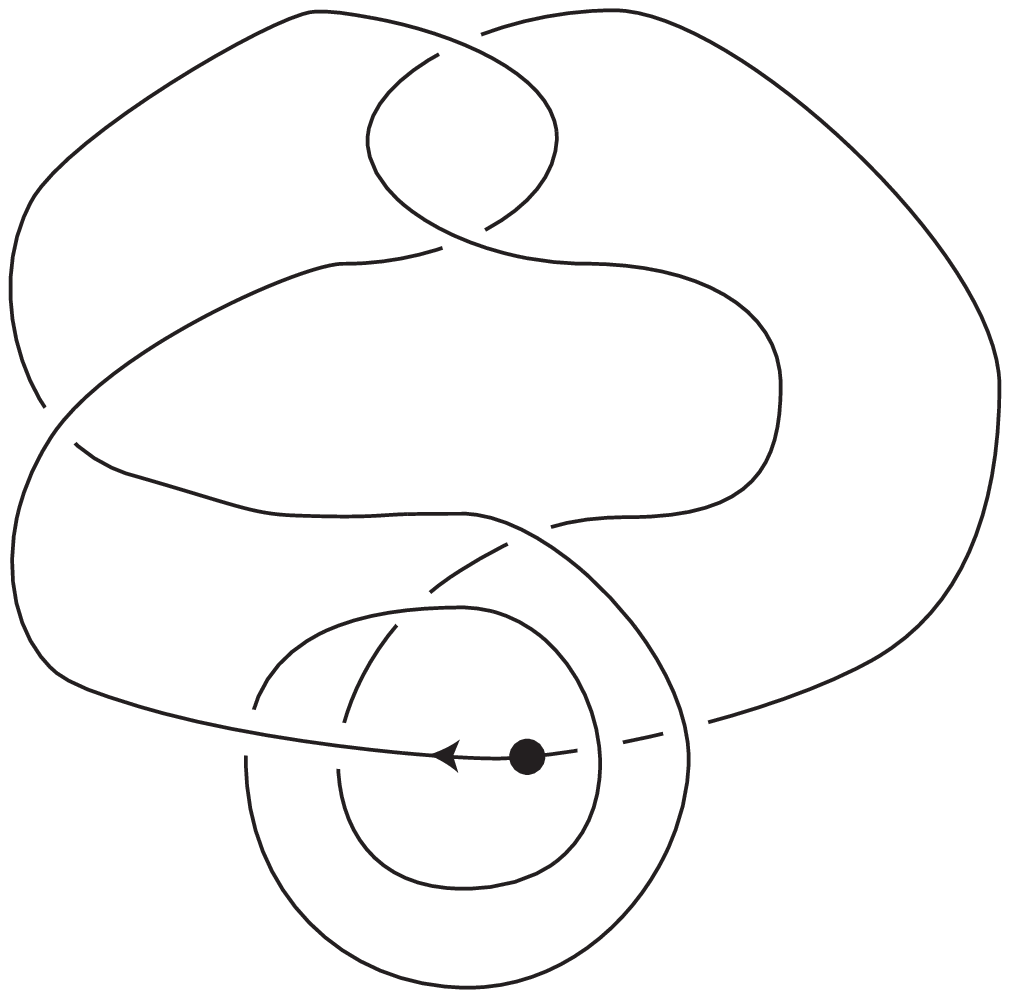}\\
	$5_2$ & $6_1$ & $7_2$\\
	\end{tabular}
	\end{center}
	\label{minimal}
\end{figure}

\begin{figure}[htbp]
	\begin{center}
	\begin{tabular}{ccc}
	\includegraphics[trim=0mm 0mm 0mm 0mm, width=.32\linewidth]{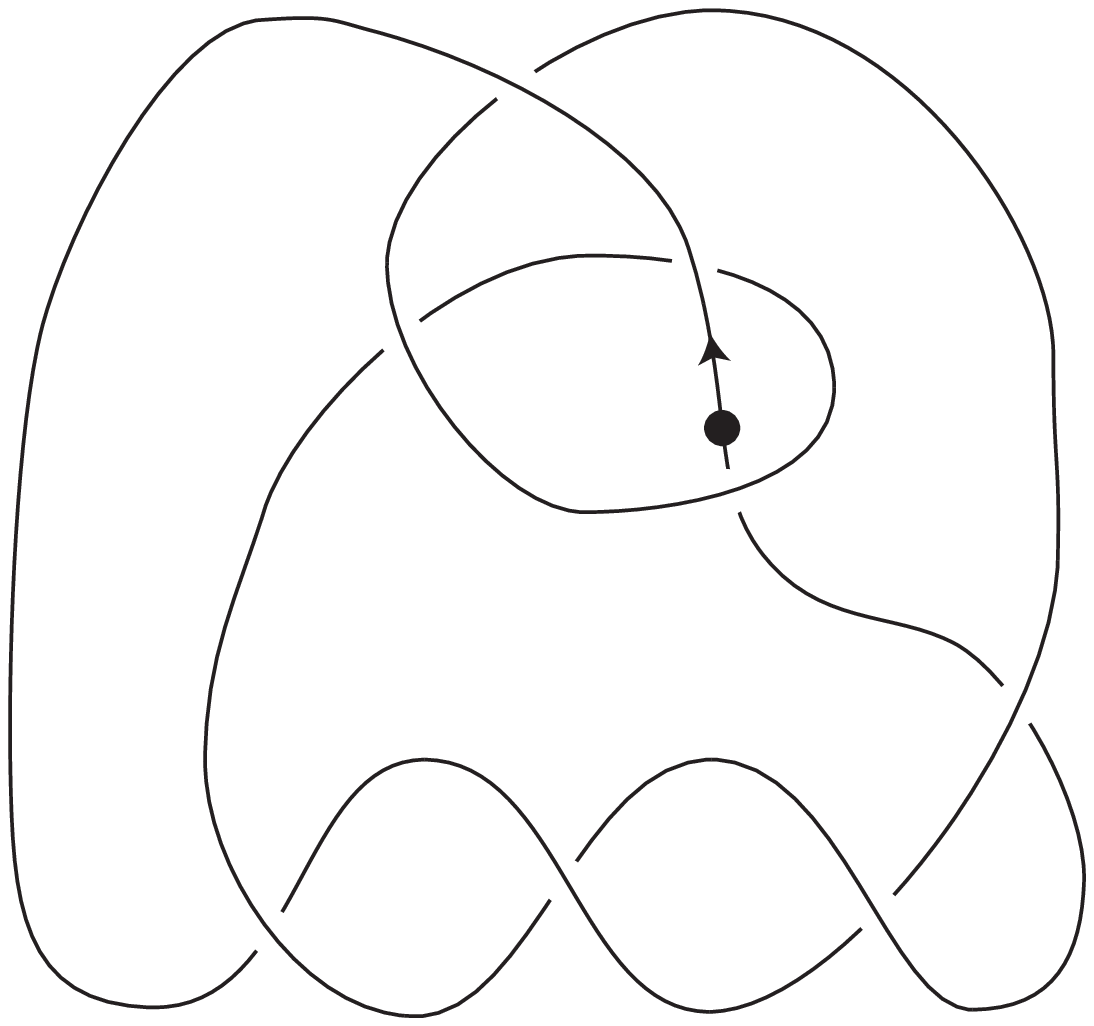}&
	\includegraphics[trim=0mm 0mm 0mm 0mm, width=.28\linewidth]{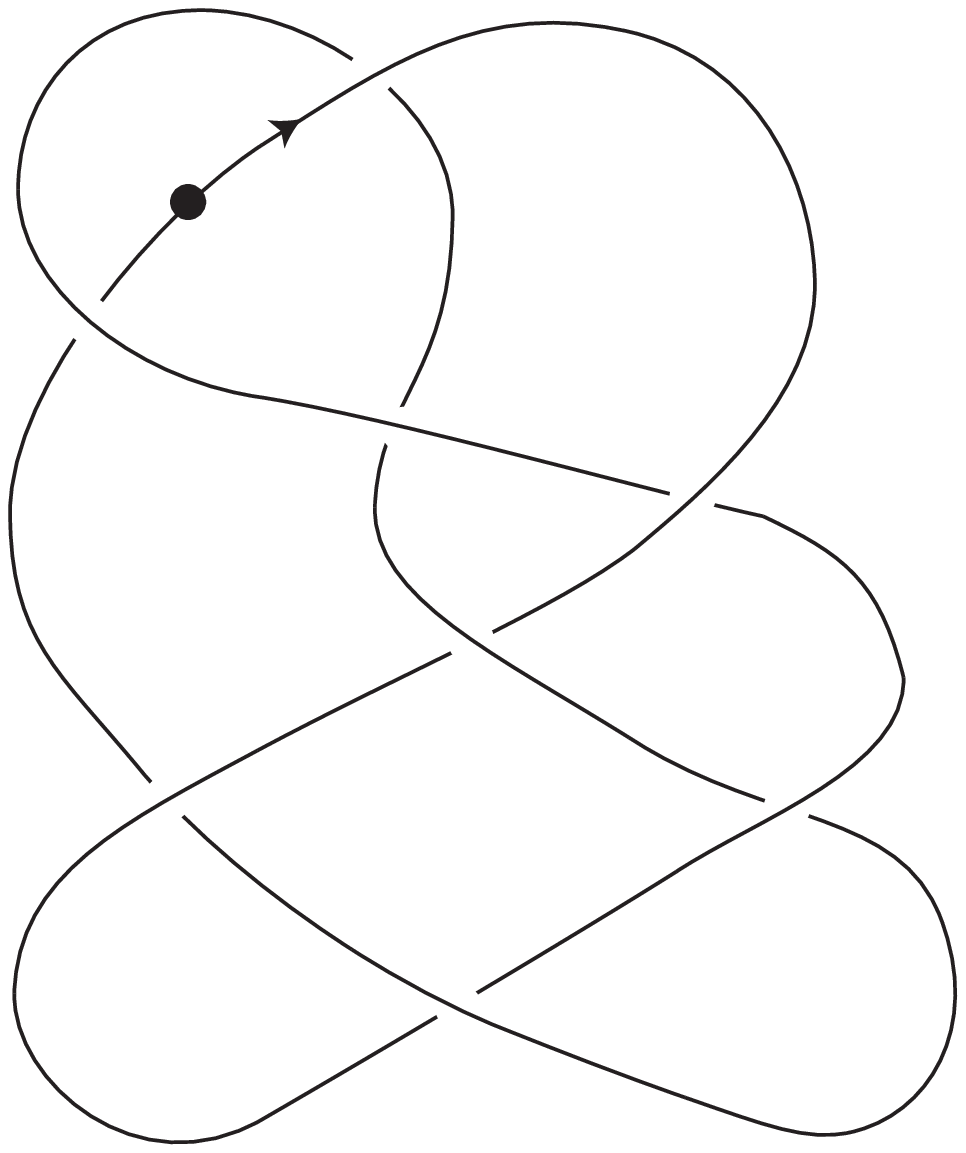}&
	\includegraphics[trim=0mm 0mm 0mm 0mm, width=.28\linewidth]{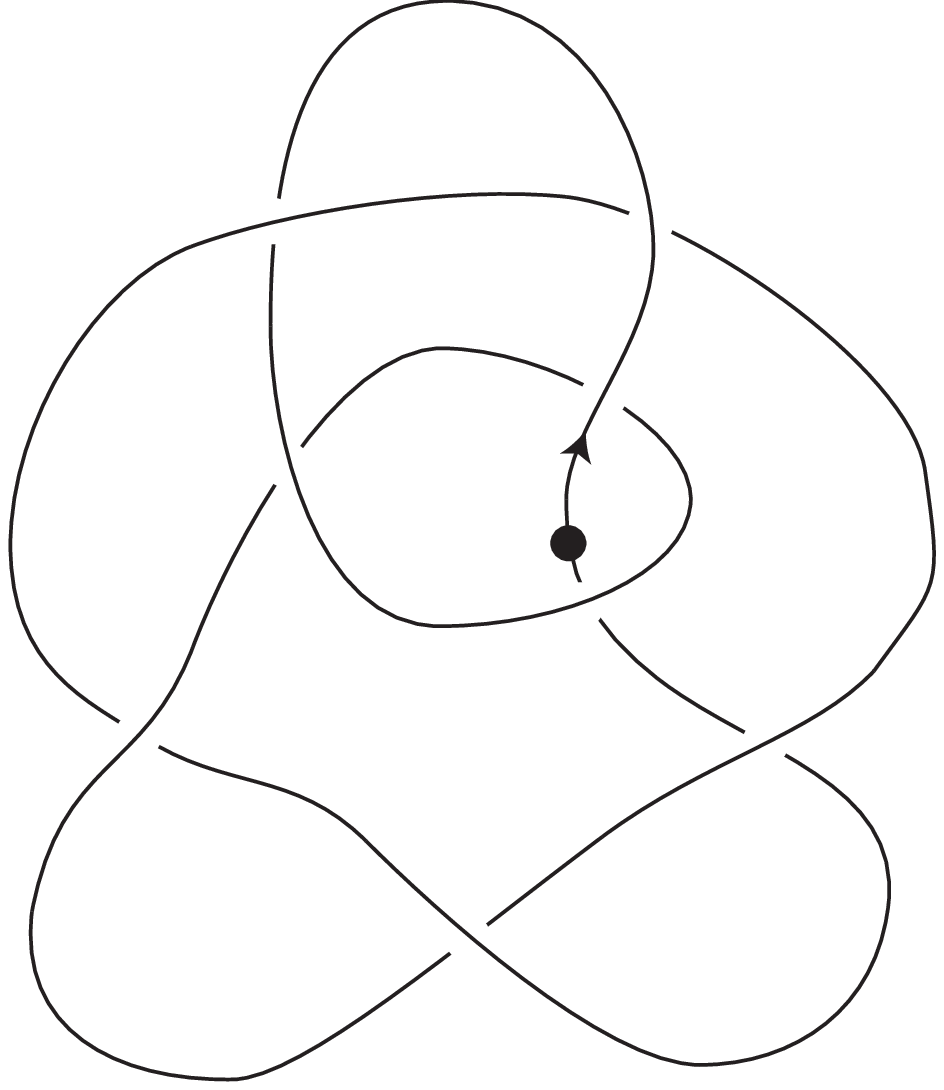}\\
	$7_3$ & $7_4$ & $7_5$\\
	\end{tabular}
	\end{center}
	\label{minimal}
\end{figure}

\begin{figure}[htbp]
	\begin{center}
	\begin{tabular}{ccc}
	\includegraphics[trim=0mm 0mm 0mm 0mm, width=.3\linewidth]{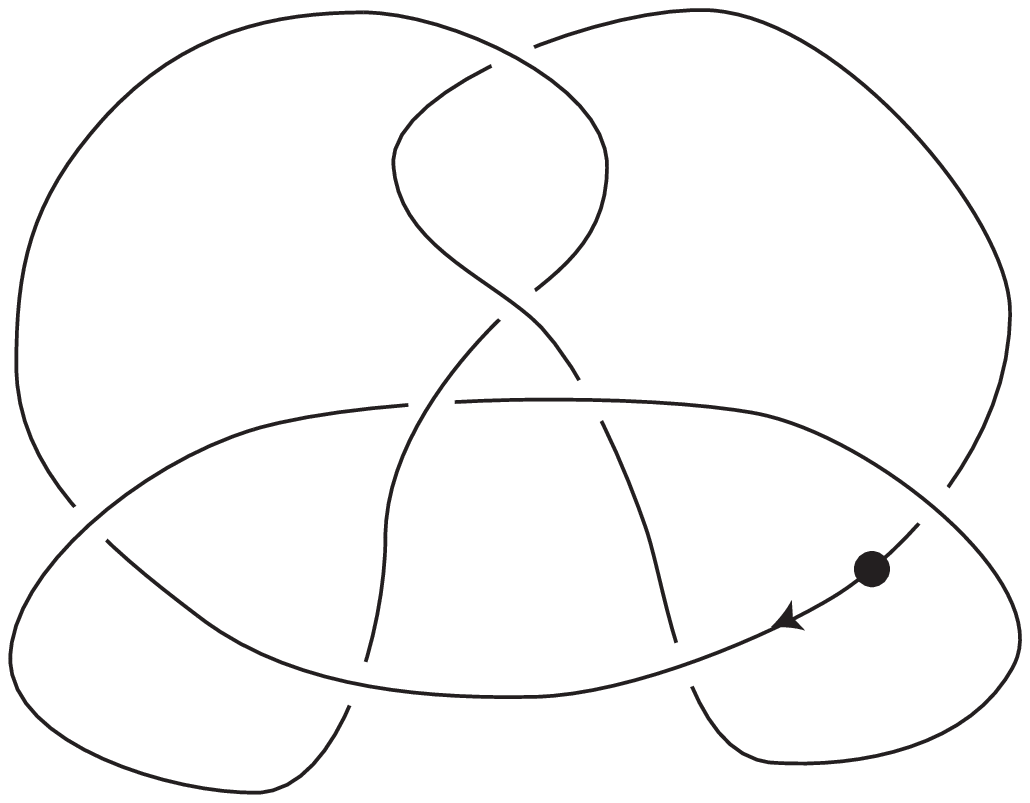}&
	\includegraphics[trim=0mm 0mm 0mm 0mm, width=.3\linewidth]{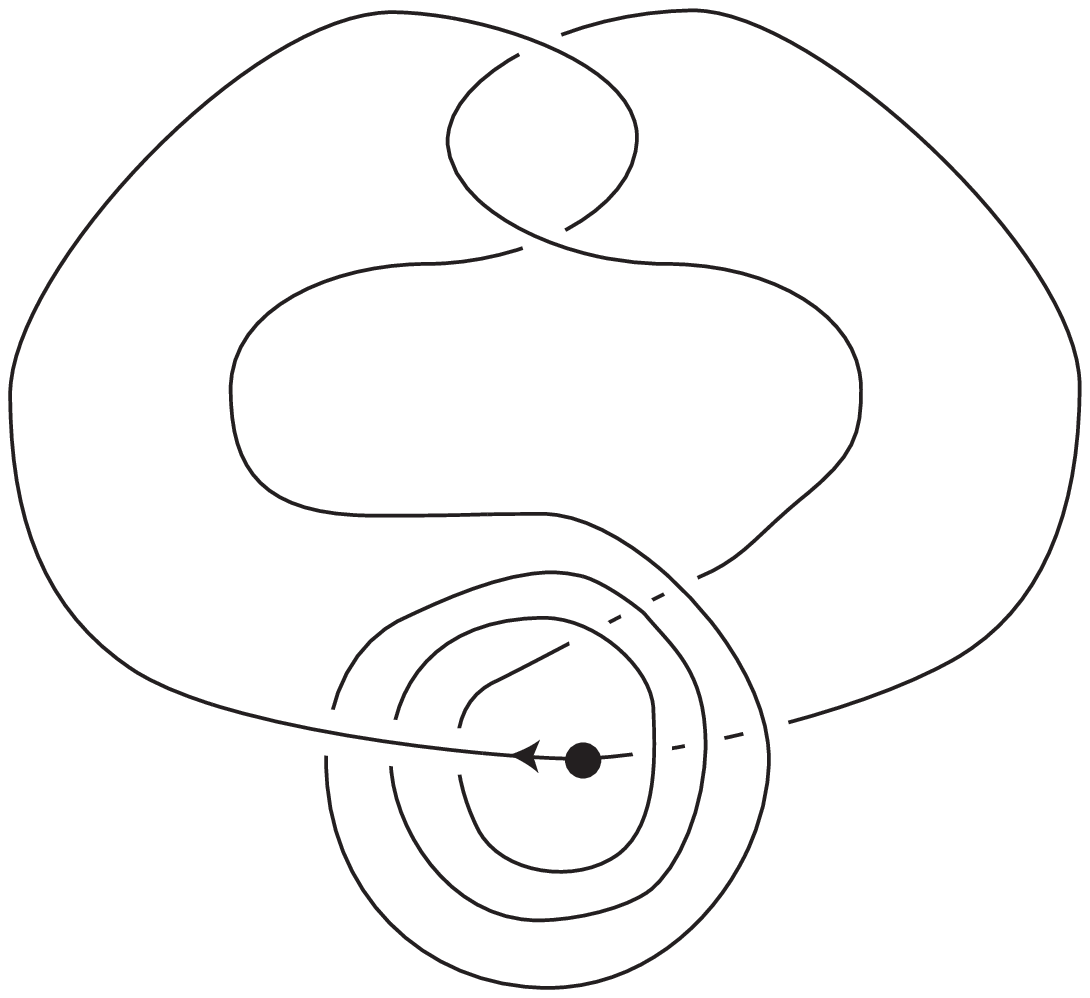}&
	\includegraphics[trim=0mm 0mm 0mm 0mm, width=.28\linewidth]{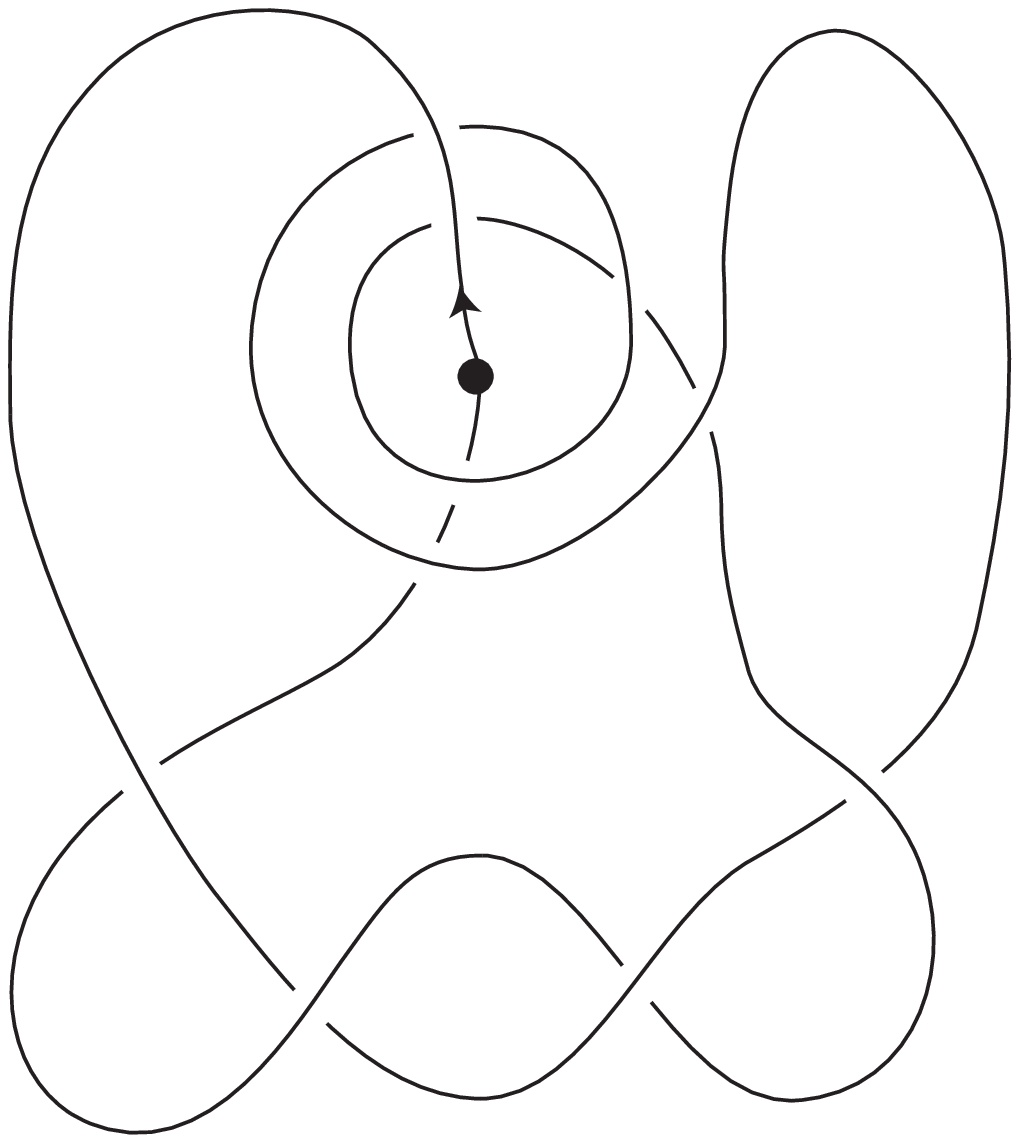}\\
	$7_6$ & $8_1$ & $8_3$\\
	\end{tabular}
	\end{center}
	\label{minimal}
\end{figure}

\begin{figure}[htbp]
	\begin{center}
	\begin{tabular}{ccc}
	\includegraphics[trim=0mm 0mm 0mm 0mm, width=.33\linewidth]{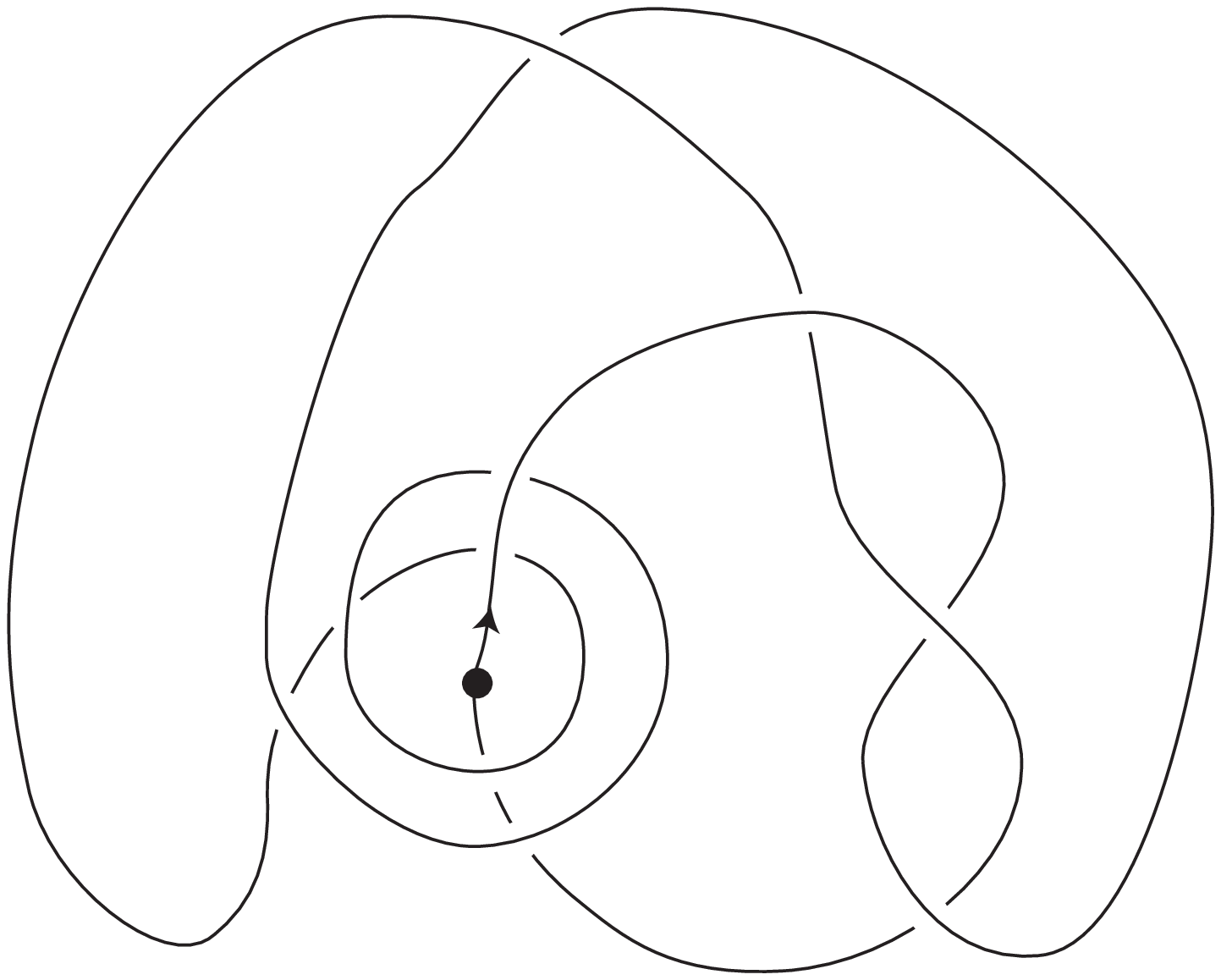}&
	\includegraphics[trim=0mm 0mm 0mm 0mm, width=.26\linewidth]{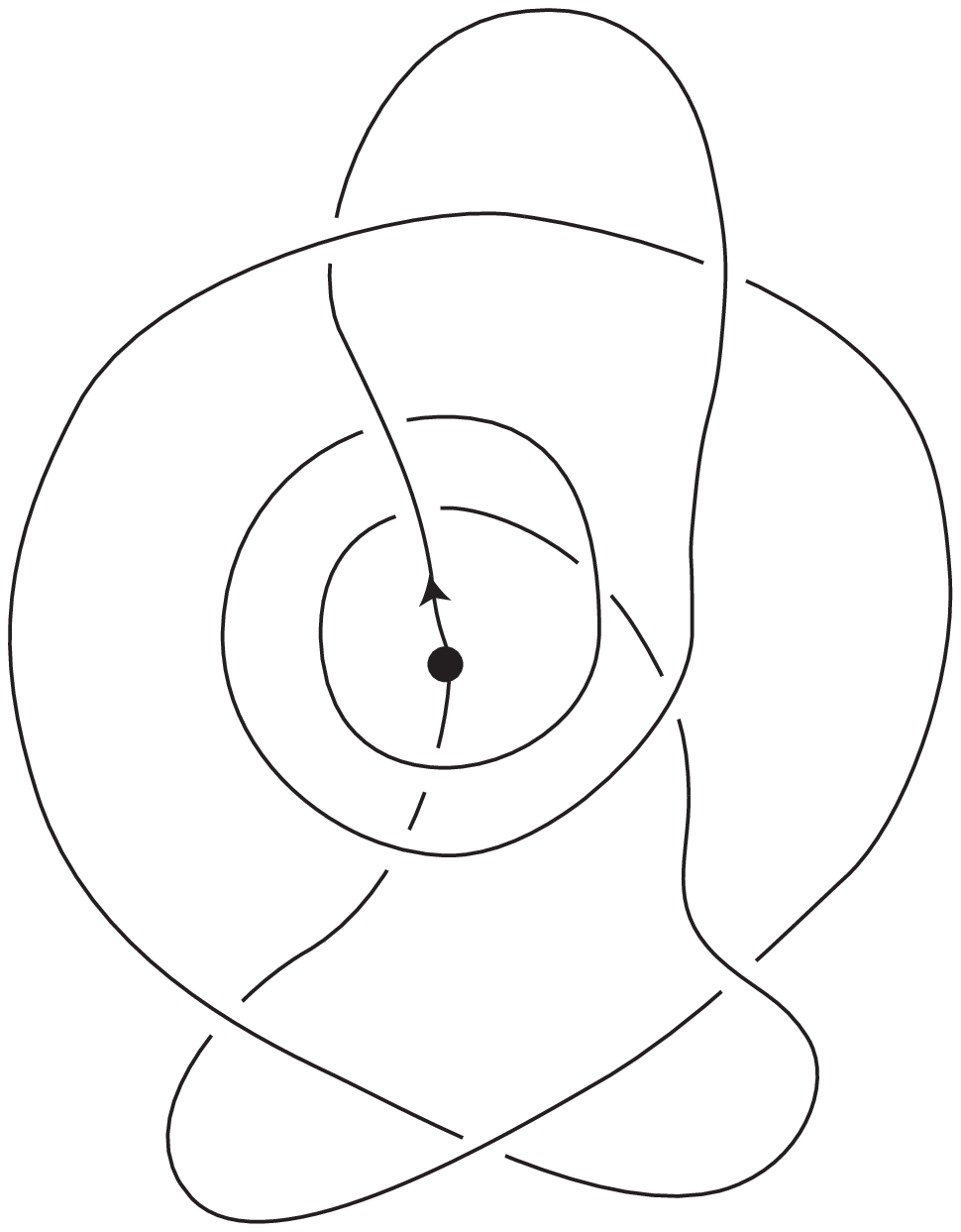}&
	\includegraphics[trim=0mm 0mm 0mm 0mm, width=.28\linewidth]{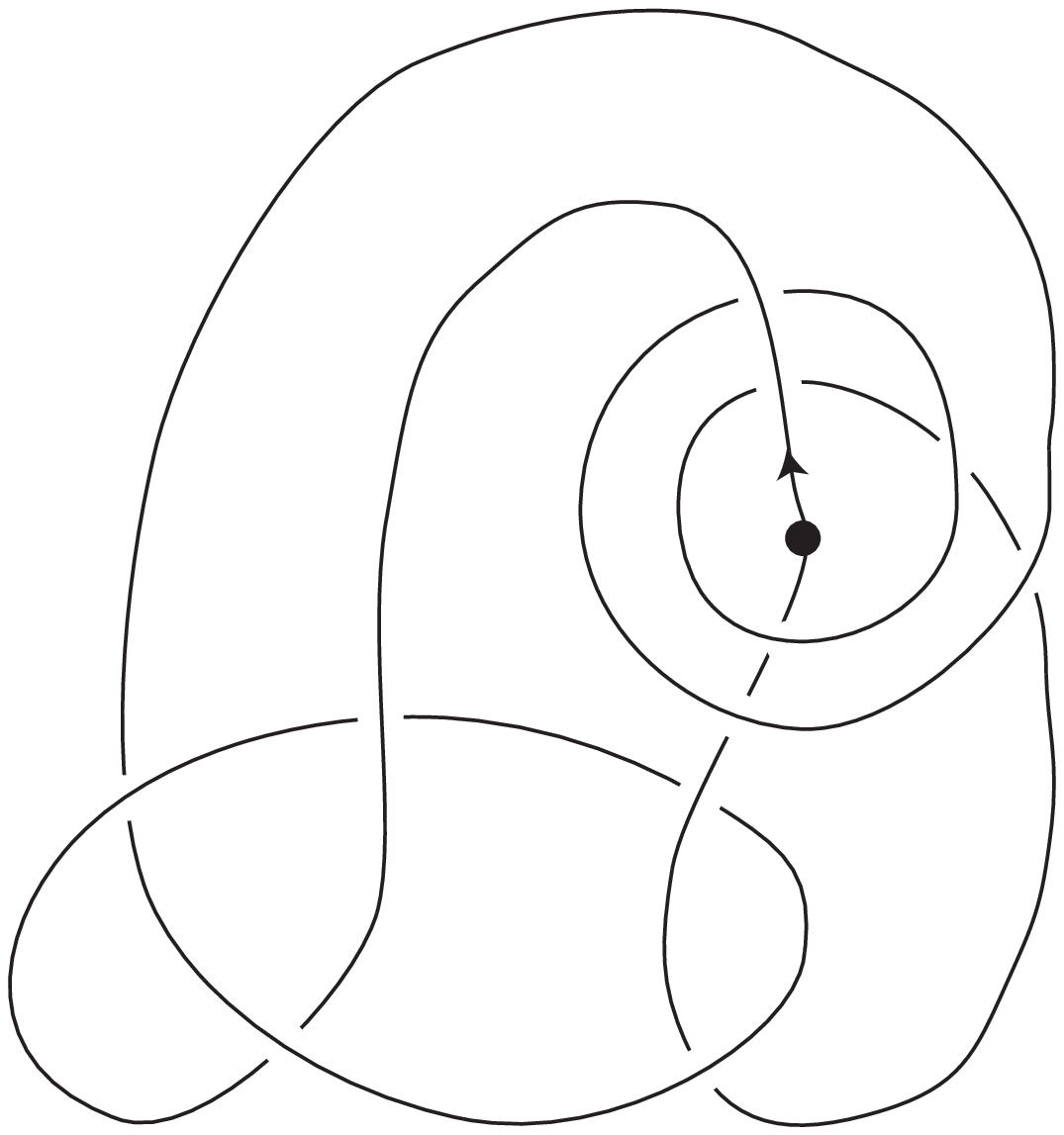}\\
	$8_4$ & $8_6$ & $8_8$\\
	\end{tabular}
	\end{center}
	\label{minimal}
\end{figure}

\begin{figure}[htbp]
	\begin{center}
	\begin{tabular}{ccc}
	\includegraphics[trim=0mm 0mm 0mm 0mm, width=.33\linewidth]{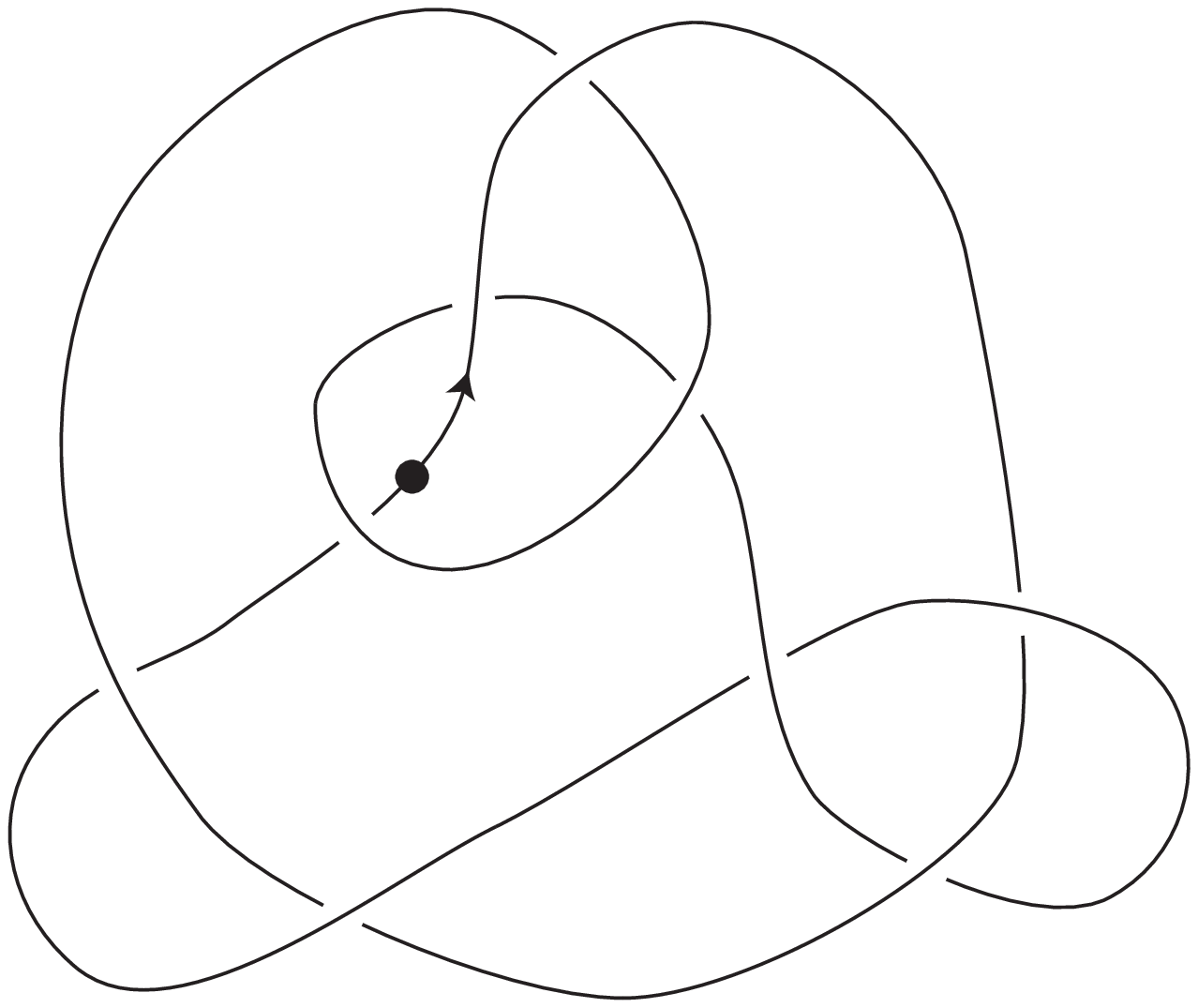}&
	\includegraphics[trim=0mm 0mm 0mm 0mm, width=.25\linewidth]{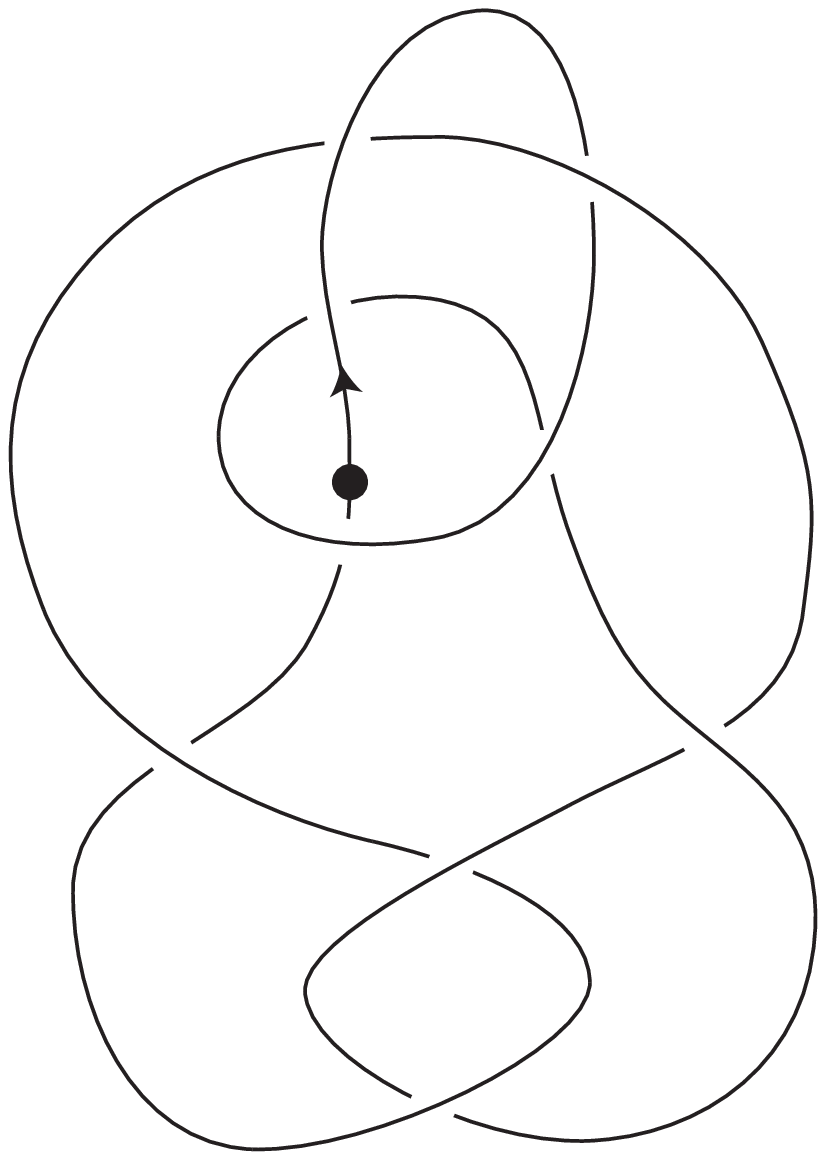}&
	\includegraphics[trim=0mm 0mm 0mm 0mm, width=.25\linewidth]{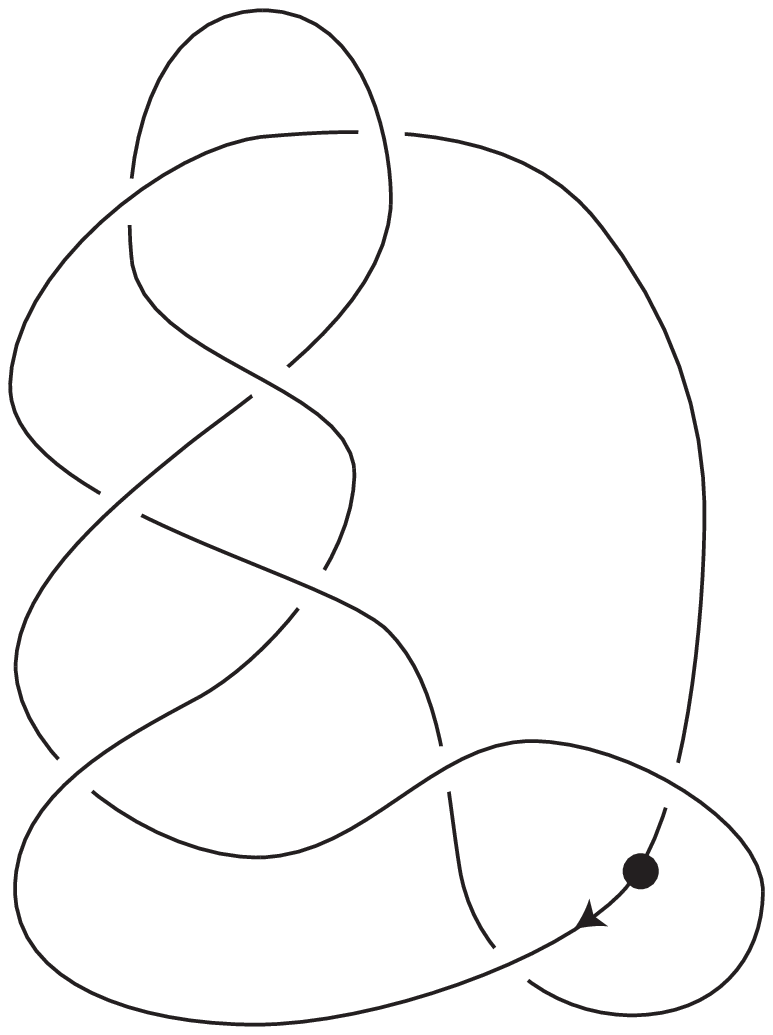}\\
	$8_{11}$ & $8_{12}$ & $8_{13}$\\
	\end{tabular}
	\end{center}
	\label{minimal}
\end{figure}

\begin{figure}[htbp]
	\begin{center}
	\begin{tabular}{ccc}
	\includegraphics[trim=0mm 0mm 0mm 0mm, width=.3\linewidth]{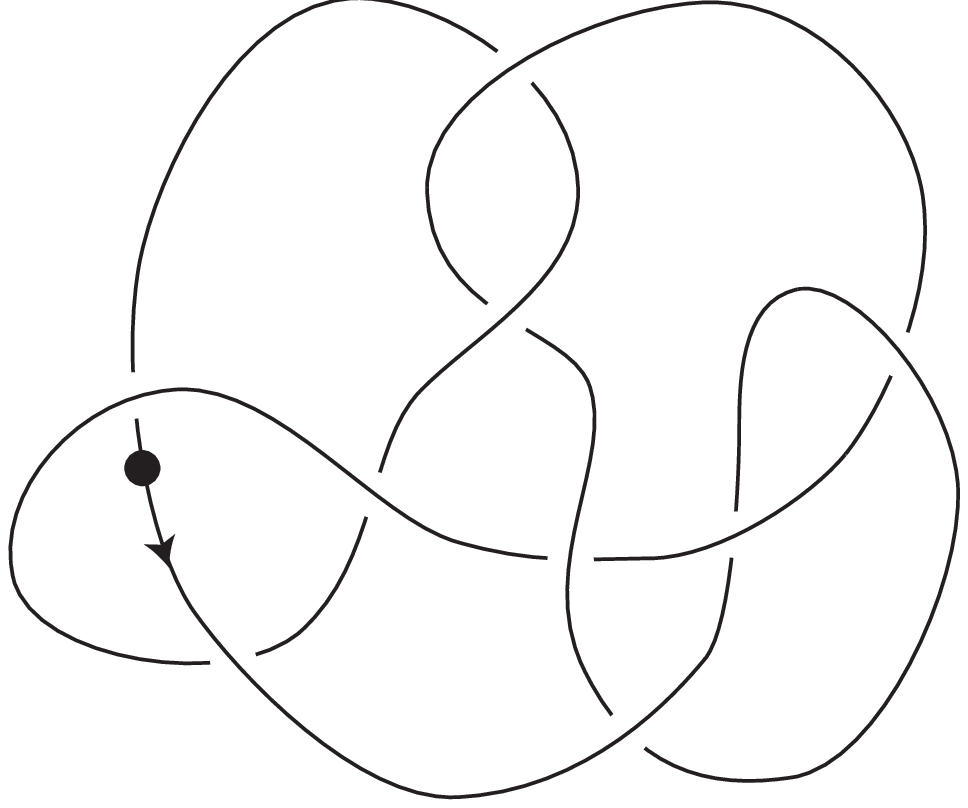}&
	\includegraphics[trim=0mm 0mm 0mm 0mm, width=.3\linewidth]{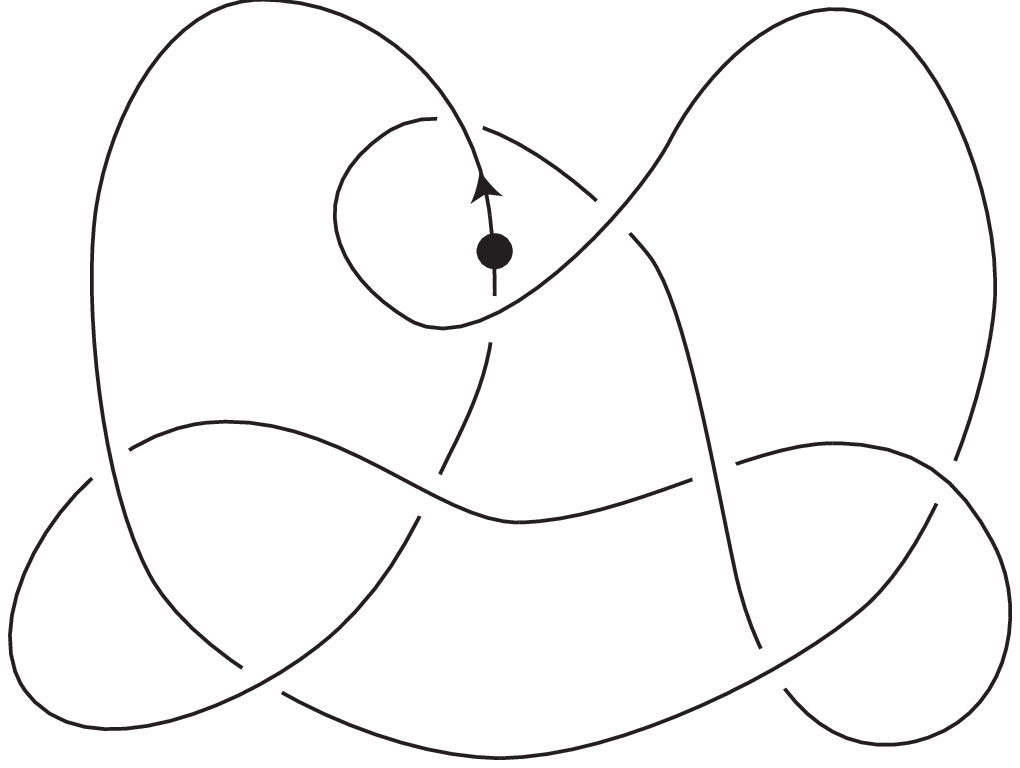}&
	\\
	$8_{14}$ & $8_{15}$ & \\
	\end{tabular}
	\end{center}
	\label{minimal}
\end{figure}

\bigskip

\noindent{\bf Acknowledgement.}
I would like to thank Ryo Nikkuni for informimg me a result of Tat Sang Fung about Theorem \ref{one}.

\bibliographystyle{amsplain}

\begin{thebibliography}{10}

\bibitem{A} C. Adams, {\em The Knot Book}, American Mathematical Society, 2004.

\bibitem{BZ} G. Burde and H. Zieschang, {\em Knots}, Walter de Gruyter, 2002.

\bibitem{C} P. Cromwell, {\em Knots and Links}, Cambridge University Press, 2004.

\bibitem{F}S. Fujimura, {\em On the ascending number of knots}, a thesis, Hiroshima University, 1998.

\bibitem{Fun} T. S. Fung, {\em Immersions in knot theory}, a dissertation, Columbia University, 1996.

\bibitem{Kaw}A. Kawauchi, {\em Survey on Knot Theory}, Birkh\"{a}user Verlag, 1996.

\bibitem{K}R. Kirby, {\em Problems in low-dimensional topology}, Geometric Topology {\bf 2}, Amer. Math. Soc., (1997) 35-473.

\bibitem{KM} P. B. Kronheimer and T. S. Mrowka, {\em Gauge theory for embedded surfaces}, Topology {\bf 32} (1993) 773-826.

\bibitem{Lic} W. B. R. Lickorish, {\em The unknotting number of a classical knot}, Contemp. Math. {\bf 44} (1985) 117-121.

\bibitem{Lic2} W. B. R. Lickorish, {\em An Introduction to Knot Theory}, Springer, 1997.

\bibitem{L} C. Livingston, {\em Knot Theory}, The Mathematical Association of America, 1996.


\bibitem{M} K. Murasugi, {\em Knot Theory and Its Applications}, translated by B. Kurpita, Birkh\"{a}user Boston, 1996.

\bibitem{O} M. Okuda, {\em A determination of the ascending number of some knots}, a thesis, Hiroshima University, 1998.

\bibitem{MO}M. Ozawa, {\em Ascending number of knots and links}, a talk in a conference ``Knot Theory'' at Waseda University, 1995.

\bibitem{OS}P. Ozsvath and Z. Szabo, {\em Knots with unknotting number one and Heegaard Floer homology}, Topology {\bf 44} (2005) 705-745.

\bibitem{R} D. Rolfsen, {\em Knots and Links}, AMS Chelsea Publishing, 2003.


\bibitem{Sch} H. Schubert, {\em \"{U}ber eine numerische Knoteninvariante}, Math. Z. {\bf 61} (1954) 245-288.


\bibitem{S} S. Suzuki, {\em An introduction to knot theory}, Saiensu-sha Co., Ltd. Publishers, 1991.


\end{thebibliography}

\end{document}